\documentclass{bioinfo}
\copyrightyear{2014}
\pubyear{2014}
\usepackage{nccbbb}

\newcommand{\mean}{\mathsf{E}}
\newcommand{\var}{\mathsf{D}}
\begin{document}
\firstpage{1}

\title[Application of the FSD for approximation of gene expression]{Application of the fractional stable distributions for approximation of gene expression profiles}
\author[V. V. Saenko, Yu. V. Saenko]{Viacheslav V. Saenko\,$^{1}$\footnote{to whom correspondence should be addressed}, Yurij V. Saenko\,$^{1}$}
\address{$^{1}$Ulyanovsk State University, Leo Tolstoy str., 42, Ulyanovsk, Russia, 432000}

\history{Received on XXXXX; revised on XXXXX; accepted on XXXXX}

\editor{Associate Editor: XXXXXXX}

\maketitle

\begin{abstract}

\section{Motivation:}
At the present time reliably established that probability density functions of gene expression of microarray experiments possess a number of universal properties. First of all these distributions have power asymptotic and secondly the shape of these distributions are inherent for all organisms and tissues. This fact led to appearance of a number works where authors are investigating various probability distributions for approximation of empirical distributions of gene expression. Nevertheless all these distributions aren't limit distribution and aren't solution of any equations. These facts by our opinion are essential shortcoming of these probability laws. Besides, expression of individual gene aren't accidental event and it depends from expression other genes. This allows to talk about existence of genic regulatory net in the cell.
\section{Results:}
In the work the class of fractional stable distributions (FSD) are described. This class of distributions is limit distribution of sum independent identical distributed random variables. These distributions have power-law asymptotic and this fact allow us to apply their for approximation of experimental densities gene expression of microarray experiments. The parameters of FSDs are statistically estimated by experimental dates and empirical density is compared whit theoretical density. In the work the algorithms of parameters estimation and simulating of FSD variables are presented. The results of such comparison allow to make conclusion that empirical densities of gene expression can be approximate by FSD.

\section{Contact:} \href{saenkovv@gmail.com}{saenkovv@gmail.com}
\end{abstract}

\section{Introduction}

The technology of hybridization DNA microarrays of high-density has opened possibility to study the expression of genes of all known genes in a single experiment. Studying the dynamics of the gene expression is one of priority trends in modern biology and medicine as it allows to understand the mechanisms of appearance of pathological conditions at the  cellular level. This means that knowledge of theoretical distribution opens outlooks in development of models of gene expression dynamics. Changing of gene expression is a complex coordinated process which depends from large number external and internal factors \cite{Macneil2011}. Therefore the probability methods most suitable for description of such processes.

Currently don't exists fixed opinion about the kind of probability distribution which describe the profiles of gene expression of microarray experiments.
Reliably established that empirical distributions are one-sided distributions, they have power-law asymptotic and character of these distributions don't changes for wide area of tissues and organisms from E. coli to H. sapiens \cite{Ueda2004a}. Such universality suggests fundamental nature of processes which leads to observable distribution of gene expression. Analogues conclusions have been obtained in works other authors \cite{LIEBOVITCH2006,Lu2009,Furusawa2003,Hoyle2002,Kuznetsov2002} where gene expression of various organisms is also investigated.

Power-law asymptotic of an experimental distribution means that theoretical distribution must have the asymptotic of following form
\begin{equation}\label{eq:Zipf}
p(x)\propto x^{-\alpha-1}, x\to\infty.
\end{equation}
In the above work  \cite{Ueda2004a} the same distribution was applied for approximation of profiles of gene expression various organisms under consideration and was showed that the parameter $\alpha$ is varying within limits from 0.69 to 1.09. In another work \cite{Furusawa2003} have investigated more than 40 tissues for 6 organisms, and for all samples the power-law distribution was obtained. In the article \cite{Kuznetsov2002} was marked that the best approximation among Poisson distribution, exponential distribution, logarithmic distribution, power-law distribution, parettolike distributions and mixtures of logarithmic and exponential distributions gives discrete Paretto distribution $p(m)=(m+b)^{-\alpha-1}/z$, where $\alpha$ is varying within limits from 0.974 to 1.88.

However the distribution (\ref{eq:Zipf}), which is named Zipf-Pareto distribution, aren't the only distribution with power-law asymptotic. In the paper \cite{Hoyle2002}  was obtained that if to make logarithmic transformation of  raw expression data and then align  and standardize them  $\xi=(\log s-\mu)/\sigma^2$, then distribution of transformed data, is well described by log-normal distribution. Here $\mu$ is mathematical expectation and $\sigma^2$ is variance of random variable $\log s$, $s$ is raw value of gene expression. In another article \cite{Lu2009} authors suggest to use double Pareto-log-normal distribution. Besides Pareto-log-normal distribution the authors in the work tested  such distribution as Zipf-Pareto distribution, log-normal distribution, log-gamma distribution, log logistic distribution, right-side Pareto-distribution. As a result, in the paper, the authors conclude that the best results are obtained with the double Pareto-lognormal distribution.

In present work the fractional stable distributions (FSD) \cite{Bening2006} are used for approximation gene expression profiles. The FSD belong to the class of infinitely divisible distributions and, in addition, they are the limit distributions of sums of independent identically distributed random variables.

Let's explain what is meant by words the limit distribution. Let we have sequence of series of independent (in each series) identical distributed random variables $\{X_{ij}, j=1,2,\dots,n_i, i=1,2,\dots\}$. Here $i$ is the number of series and each series contain $n_i$ of random variables.  Form the sums
\begin{equation}\label{eq:sumSeries}
S_i=\sum_{j=1}^{n_i}X_{ij},\quad i=1,2,\dots.
\end{equation}
Of course, in common case, random variables $X_{ij}$ must be aligned and standardized, but consideration of this question goes beyond the bounds of the article. We assume that the random  variables $X_{ij}$ have been normalized in a certain way and centered. Then the distribution $p(x)$ will be the limit distribution if it will be distribution of the sum $S_i$. Depending on imposed on the random variables $ X_{ij} $ properties result in different classes of limit distributions.

Widely known limiting distribution is a normal distribution
\begin{equation}\label{eq:normalLaw}
p(x)=(4\pi\sigma^2)^{-1/2}\exp\left(-(x-\mu)^2/(4\sigma^2)\right),
\end{equation}
where $\mu=\mean X_{ij}$ is mathematical expectation and $\sigma^2=\var X_{ij}$ is variance of random variable $X_{ij}$. According central limit theorem in order to obtain the normal distribution as the limit for a series of sums (\ref{eq:sumSeries}) sufficient to require that the random variables $ X_ {ij} $ have a finite mathematical expectation  and finite variance $\mean X_{ij}<\infty$, $\var X_{ij}<\infty$. If we now require that the variance of random variables $X_{ij}$ is infinite ($\var X_{ij}=\infty$) and distribution of each variable has asymptotic (\ref{eq:Zipf}) then the limit distribution of the sum (\ref{eq:sumSeries}) will belong to the class of stable laws \cite{Zolotarev1986, Uchaikin1999}. It should be noted that another necessary condition of appearance of the class stable distributions as the limit distributions is that number of terms in the sum   (\ref{eq:sumSeries}) were fixed but enough large. If suppose that the number of terms in each sum (\ref{eq:sumSeries}) are random then the limit distribution of sum $S_i, i=1,2,\dots$ will be belong to the class FSD.

The classes of stable and FSDs are well-known but rarely used. The reason for this lies in the lack of expressions for the probability density function presented in terms of elementary functions. Nevertheless, it exist numerical methods of calculation of value of probability density function in a given point \cite{Uchaikin2002} and methods of statistical estimation of parameters of FSD \cite{Bening2004}. In the present work profiles of gene expression are approximated by FSD on the basis of these methods. Parameters of distributions are statistically estimated according to raw data and then theoretical densities are compared with empirical densities.

\section{Fractional stable distributions}

In first time FSD was introduced in work \cite{Kotulski1995} where the scheme of random walks with traps has been considered. Let in initial time moment $t=0$ particle appears and makes instantaneous jump on distance $X_1$. After that rests in this position for random time $T_1$. After that a particle again makes instantaneous jump on distance $X_2$ after that all process repeats same way. Random variables $T_1, T_2,\dots$ and $X_1,X_2,\dots$ are independent identically distributed and these two sequences aren't depend from each other. As result we obtain the particle trajectory. Once the particle ends his story, next particle trajectory is constructed.

As results a coordinate of each particle for $i$-th trajectory can be describes by the sum
\begin{equation}\label{eq:sumFSD}
S_i=\sum_{j=1}^{N_i(t)}X_{ji},\quad i=1,2,\dots,
\end{equation}
where the number of term $N_i(t)$ in the sum are defined from the process
\begin{equation}\label{eq:sumTime}
\sum_{j=1}^{N_i(t)}T_{ji}<t\leqslant\sum_{j=1}^{N_i(t)+1}T_{ji},\quad i=1,2,\dots.
\end{equation}
If now to suppose that distribution $\rho_X(x)$ of random variables $X_j$ and distribution $\eta_T(t)$ of random variables $T_j$ have power-law asymptotic of form $\rho_X(x)\propto \alpha x_0^\alpha x^{-\alpha-1}, 0<\alpha\leqslant2, x\to\infty$, $\eta_T(t)\propto \beta t_0^\beta t^{-\beta-1}, 0<\beta\leqslant1, t\to\infty$, then the limit distribution of the sum (\ref{eq:sumFSD}) are described by expression
$$
p(x,t)= t^{-\beta/\alpha} q(xt^{-\beta/\alpha};\alpha,\beta,\theta,\lambda),
$$
where $q(x;\alpha,\beta,\theta,\lambda)$ is FSD \cite{Kolokoltsov2001,Bening2006}. FSD are expressed through Mellin's transform of two stable distributions
\begin{equation}\label{eq:FSD}
q(x;\alpha,\beta,\theta,\lambda)=\int g(xy^{\beta/\alpha};\alpha,\theta,\lambda)g(y;\beta,1,1)y^{\beta/\alpha}dy.
\end{equation}
Here $g(x;\alpha,\theta,\lambda)$ is density function of strictly stable law and $g(y;\beta,1,1)$ is density of one-sided strictly stable law with characteristic function \cite{Zolotarev1986}
\begin{equation}\label{eq:formc}
\hat g(k;\alpha,\theta;\lambda)=\exp\left(-\lambda |k|^\alpha\exp\left(-i\alpha\theta(\pi/2)\mathrm{sign} k\right)\right).
\end{equation}
The strictly stable distributions are fully defined by their three parameters $\alpha,\theta,\lambda$, where  $\alpha\in(0,2]$ is characteristic parameter, $\theta$ is asymmetry parameter $(|\theta|\leqslant\min(1,2/\alpha-1))$ and $\lambda>0$ is scale parameter.

As we can see from (\ref{eq:FSD}) FSD are defines by four parameters. The parameters $\alpha$ and $\beta$ are two characteristic parameters. Parameters $\alpha$ and $\beta$ are varying  in the limits $\alpha\in(0,2]$ $\beta\in (0,1]$. Domain of variation of parameters $\theta$ and $\lambda$ coincides with domain of variation respective parameters of stable distribution and they have the same meaning. The FSD has a power-law asymptotic (\ref{eq:Zipf}). When $\beta=1$ class of FSD pass to the class strictly stable distributions. Indeed when $\beta=1$ and $\theta=1$ strictly-stable law $g(y;1,1)$ is singular distribution at point $y=1$.  Hence, from (\ref{eq:FSD}) we obtain
$
\int_0^\infty g(xy^{\beta/\alpha};\alpha,\theta,\lambda)\delta(y-1)y^{\beta/\alpha}dy=g(x;\alpha,\theta,\lambda).
$
In the case when $\alpha=2$, $\beta=1$ è $\theta=0$ from (\ref{eq:FSD}) and (\ref{eq:formc}) we obtain that FSD is passes to the normal distribution (\ref{eq:normalLaw}). Hence, the class of fractional stable laws involves the class of stable distributions  to
which include the Gaussian distribution, the Cauchy distribution,
Levy-Smirnov distribution.

The fact that the FSD have power asymptotic behavior suggests the possibility of using this class of distributions to describe the density distribution of gene expression levels. Therefore make the following assumption. Suppose, that gene expression levels obtained with each probe experiments with chip microarrays represent a sample of random variables $ Z_i, i = 1 , \dots, N $ with fractional stable distribution $ q (x; \alpha, \beta, \theta , \lambda) $. Therefore, to approximate the experimental distribution of expression levels of genes it is necessary for the sample $ Z_i $ evaluate parameters $ \alpha, \beta, \theta, \lambda $.

\section{Estimation of FSD parameters}\label{sec:FSDEst}

There is only two method of estimation of parameters of FSD in present time. The first method has been described \cite{Saenko2012} and estimator is obtained on basis maximum likelihood method. The second method has been described in the work \cite{Bening2004} and estimators were obtained on basis the moment method.  However both of the methods were found to be not convenient for estimation of the parameters of FSD of gene expression of microarrays. Indeed, it is necessary to know analytical expression for FSD for calculation of likelihood function at usage the first method. In the article \cite{Saenko2012} this difficulty succeeded avoiding by usage of local estimator of Monte Carlo method for calculation of symmetric FSD. However this estimator can be used for estimation of the parameters only of symmetric FSD while gene expression distributions are one-sided distribution. Usage of the second estimator  has shown that the parameters are not estimated correctly. For this reason another estimator of the parameters of FSD has been developed.

The basis of this algorithm the idea of minimizing of the distance $d(\mathbf{P}_\mathbf{\Theta}, \mathbf{Q})$ between two distributions  is underlain. Here $\mathbf{Q}$ is empirical distribution the parameters of which necessary to estimate and $\mathbf{P}_\mathbf{\Theta}$ is theoretical distribution the parameters $\mathbf{\Theta}$ known.

As known \cite{Kolokoltsov2001} the fractional stable random variable can be represented as ratio of two strictly stable random variable
\begin{equation}\label{eq:fsrv}
Z(\alpha,\beta,\theta,\lambda)=\lambda Y(\alpha,\theta)[S(\beta,1)]^{-\beta/\alpha},
\end{equation}
where $Y(\alpha,\theta)$ -- stable random variable and $S(\beta,1)$ -- one-sided stable random variable are distributed according to stable law with characteristic function (\ref{eq:formc}). Since the algorithms of simulation of stable random variables well known (see Appendix~\ref{app:FSRV}) we can simulate fractional stable random variable $Z(\alpha,\beta,\theta,\lambda)$ without any difficulty and as result one can construct a histogram of distribution. Thus the task is reduced to calculation of the distance $d(\mathbf{P}_\mathbf{\Theta}, \mathbf{Q})$ between histogram of FSD and the histogram of empirical distribution which constructed according to the sample $Z_1,Z_2,\dots,Z_N$. Here the sample $Z_i,i=1,\dots,N$ is experimental data of gene expression  obtained from microarray experiments. As a measure of such distance one can choose
\begin{equation}\label{eq:chi2dist}
d(\mathbf{P}_\mathbf{\Theta},\mathbf{Q})=
\sum_{i=1}^n\frac{(N\mathbf{P}_\mathbf{\Theta}(\Delta_i)-\nu_i)^2}
{N\mathbf{P}_\mathbf{\Theta}(\Delta_i)},
\end{equation}
where $\Delta_1,\dots,\Delta_n$ is partition of the domain under consideration $R\equiv\{x: a\leqslant x\leqslant b\}$ on $n$ disjoint intervals at the same $\cup_{i=1}^n\Delta_i=R$ and $\nu_i$ is number of observations fallen into interval $\Delta_i$. The distance (\ref{eq:chi2dist}) is called $\chi^2$ distance. As a result the estimator of $\hat{\mathbf{\Theta}}\equiv(\hat\alpha,\hat\beta,\hat\theta,\hat\lambda)$ of parameters  $\mathbf{\Theta}\equiv(\alpha,\beta,\theta,\lambda)$ will be those values of the $\mathbf{\Theta}$ at whose the distance (\ref{eq:chi2dist}) takes a minimal value.

The optimization algorithm of Hooke-Jeeves's \cite{Bunday1984} was used  for minimization of the distance (\ref{eq:chi2dist})  according to parameters $\mathbf{\Theta}$. We note here some of the features of optimization algorithms that must be considered in their application. Any optimization algorithm based on the idea of comparing values of the studied functional  in two neighboring points $\mathbf{\Theta}_i$  and $\mathbf{\Theta}_{i+1}$.  The algorithm starts from the initial point $\mathbf{\Theta}_0$ and calculates the values of the distance (\ref{eq:chi2dist}) at points $\mathbf{\Theta}_i$ and $\mathbf{\Theta}_{i+1}$, $i=0,1,\dots$.
In the case if $d(\mathbf{P}_{\mathbf{\Theta}_{i+1}},\mathbf{Q})< d(\mathbf{P}_{\mathbf{\Theta}_i},\mathbf{Q})$ (for the case when the functional is minimized) then algorithm moves to the point $\mathbf{\Theta}_{i+1}$. Thus, during optimization process the algorithm passes through points sequence $\mathbf{\Theta}_0\to\mathbf{\Theta}_1\to\dots\to\mathbf{\Theta}_m$ and in each of them the values (\ref{eq:chi2dist}) are calculated. At the same time the number of points in this sequence is random but finite.  The point $\mathbf{\Theta}_0$ is called initial point and it position is chosen arbitrarily. From this it becomes evident the nearer we take position of the point $\mathbf{\Theta}_0$ to the minimum (maximum) of the functional, the faster the algorithm will find position of the minimum (maximum). With regard to the task of estimation of the parameters of FSD then a  point  $\mathbf{\Theta}$ is described by four coordinates $(\alpha,\beta,\theta,\lambda)$ and the task consists in minimization of the distance (\ref{eq:chi2dist}) according to these four parameters. Hence we can use for determination of position of initial point  the algorithm of estimation of the parameters of FSD (see Appendix~\ref{sec:prmestmom}) was obtained in the work \cite{Bening2004}.

As was noted above during minimization process of the distance (\ref{eq:chi2dist}) the algorithm moves along a trajectory $\mathbf{\Theta}_0\to\dots\to\mathbf{\Theta}_i\to\dots\to\mathbf{\Theta}_m$ at the same each point $\mathbf{\Theta}_i$  on $i$-th step is described by four coordinates $(\alpha_i,\beta_i,\theta_i,\lambda_i)$. At the same time the parameters of FSD can take values from the domain
$G=\{(\alpha,\beta,\theta,\lambda): 0<\alpha\leqslant2,0<\beta\leqslant1,-1\leqslant\theta\leqslant1,\lambda>0\}$
Outside of the domain $G$ the FSD isn't defined. As a results, it is possible that the algorithm may go out beyond $G$. More precisely this means that one or more coordinate of point $\mathbf{\Theta_i}$ may go out beyond of $G$. In order to avoid such situation it is necessary to introduce penalty function. The main idea of penalty function consists in that, what this function increases (at minimization process) the value of the functional (\ref{eq:chi2dist}) when the algorithm goes out beyond $G$ and it doesn't change of the functional if $\mathbf{\Theta}_i\in G$. As such function the following function was chosen
$$
f(\Theta;\tilde\Theta)=\left\{\begin{array}{cl}
\exp(A|\tilde\Theta-\Theta|),&\Theta\notin G\\
1,&\Theta\in G,
\end{array}\right.
$$
where $\Theta$ denotes one of the parameters $\alpha,\beta,\theta,\lambda$, $\tilde\Theta$  denotes the closest to the $\Theta$ boundary point with respect the corresponding parameter and $A$ some multiplier ($A\gg1$). As seen the penalty is introduced with respect each of the parameters. As a results we can introduce the helper functional $\mathcal{D}(\mathbf{P}_{\mathbf{\Theta}},\mathbf{Q})$ for calculation of the distance between the theoretical distribution $\mathbf{P}_{\mathbf{\Theta}_i}$ and the empirical distribution $\mathbf{Q}$. This functional can be represented in the form
\begin{equation}\label{eq:helperFun}
\mathcal{D}(\mathbf{P}_{\mathbf{\Theta}},\mathbf{Q})=
\left\{\begin{array}{ll}
d(\mathbf{P}_{\mathbf{\Theta}},\mathbf{Q})+f(\mathbf{\Theta},\tilde{\mathbf{\Theta}})-1,& \mathbf{\Theta}\notin G,\\
d(\mathbf{P}_{\mathbf{\Theta}},\mathbf{Q}),&\mathbf{\Theta}\in G,
\end{array}\right.
\end{equation}
where $f(\mathbf{\Theta},\tilde{\mathbf{\Theta}})=f(\alpha;\tilde\alpha)f(\beta;\tilde\beta)
f(\theta;\tilde\theta)f(\lambda;\tilde\lambda).$
As seen the functional (\ref{eq:helperFun}) redefines the functional (\ref{eq:chi2dist}) to wider domain $G'=\{(\alpha',\beta',\theta',\lambda'):-\infty<\alpha'<\infty,-\infty<\beta'<\infty,
-\infty<\theta'<\infty,-\infty<\lambda'<\infty,\}$, at the same time $G\subset G'$. Besides it introduces the penalty for going out of the optimization algorithm beyond $G$. This penalty is the greater, the greater distance algorithm goes beyond $G$. It is seen from (\ref{eq:helperFun}) the value of $\mathcal{D}(\mathbf{P}_{\mathbf{\Theta}},\mathbf{Q})$ is increased in comparing of $d(\mathbf{P}_{\mathbf{\Theta}},\mathbf{Q})$ if $\mathbf{\Theta}\notin G$. Since we find minimum of functional $d(\mathbf{P}_{\mathbf{\Theta}},\mathbf{Q})$ the this forces the algorithm to return into the domain $G$.

\section{Results}

As was noted above there are two facts which allow us to make an assumption about fractional stable nature of distribution of an gene expression.
First of all the distribution of gene expression has power-law asymptotic $\propto x^{-\alpha-1}$. The FSDs have exactly the same asymptotic. Secondly, the shape of the gene expression distribution is very similar to the shape of the FSD. Consequently, the next step is testing the hypothesis about fractional stable nature of gene expression distribution. There are more fundamental reasons which may lead to power-law distributions.  The gene expression in a cell is coordinated process and large groups of genes may change their expression in dependence from expression of others genes. Most genes in a cell are grouped into special groups - signaling or metabolic pathways. At present are revealed  more than 2100 signaling and metabolic pathways. If at some time moment particular gene activates and it begins to synthesize its mRNA then activates immediately a set of genes associated with this gene. As a result a concentration of a connected set of mRNA may sharply increase and as consequence the intensity of a emission sharply grows. At the same time expression of another group of genes may be suppressed. Such variation of gene expression must lead to power-law distributions.

Since there are many manufacturers of microarrays at present time then for our aims there were selected experimental data obtained from microarrays of three manufacturers: Affymetrix, Agilent è Illumina. For microarrays of the Affymetric company were processed of gene expression the following organisms mammals (human and rat), bird (chiken), worms (C. elegans), plant (rice and Arabidopsis thaliana), insect (drosophila), bacterium (P. aeruginosa). For microarrays of the Agilent company were processed of gene expression the following organisms: mammals, fish, bird, plant, insect, bacterium and fungus. For microarrays on the Illumina company were processed three organisms: human and rat. All experimental data were obtained from free databases ArrayExpress (http://www.ebi.ac.uk/arrayexpress/) and Gene Expression Omnibus (http://www.ncbi.nlm.nih.gov/geo/).

We were interested of data which had not been exposed any primary processing (data from RAW files).
The channels PM and AM were processed separately for microarrays of the Affymetrix company.
The red and green channels were processed for the Agilent microarrays. In particular the following channels were processed: red median signal (rMedianSignal), green median signal (gMedianSignal), red mean signal (rMeanSignal), green mean signal (gMeanSignal). For microarrays of the Illumina company RAW data were processed. All expression which were processed weren't undergo any preliminary normalization or processing.

The process of processing looks as follows. Expression for the organism under consideration from processed channel is considered as sample of independent identical distributed random variable $Z_1, Z_2,\dots, Z_N$. The parameters $\hat\alpha,\hat\beta,\hat\theta,\hat\lambda$ of the FSD are estimated under this sample by algorithm which has been described in the Section~\ref{sec:FSDEst}. After that a sample of fractional stable random variables $Z(\alpha, \beta,\theta,\lambda)$ were simulated with values of the parameters $\hat\alpha,\hat\beta,\hat\theta,\hat\lambda$ had been estimated. For simulation random variables $Z(\alpha, \beta,\theta,\lambda)$ the algorithm described in Appendix~\ref{app:FSRV} was used. Next histogram was constructed. At the same time on sample $Z_1, Z_2,\dots, Z_N$ a histogram of gene expression levels was constructed. After this theoretical and empirical distributions were compared. In the case when these distributions differed insignificantly then for the $\chi^2$ Pirson's criterion was applied for checking the hypothesis about coincidence of  these two distributions.

\begin{figure}
\includegraphics[width=0.24\textwidth]{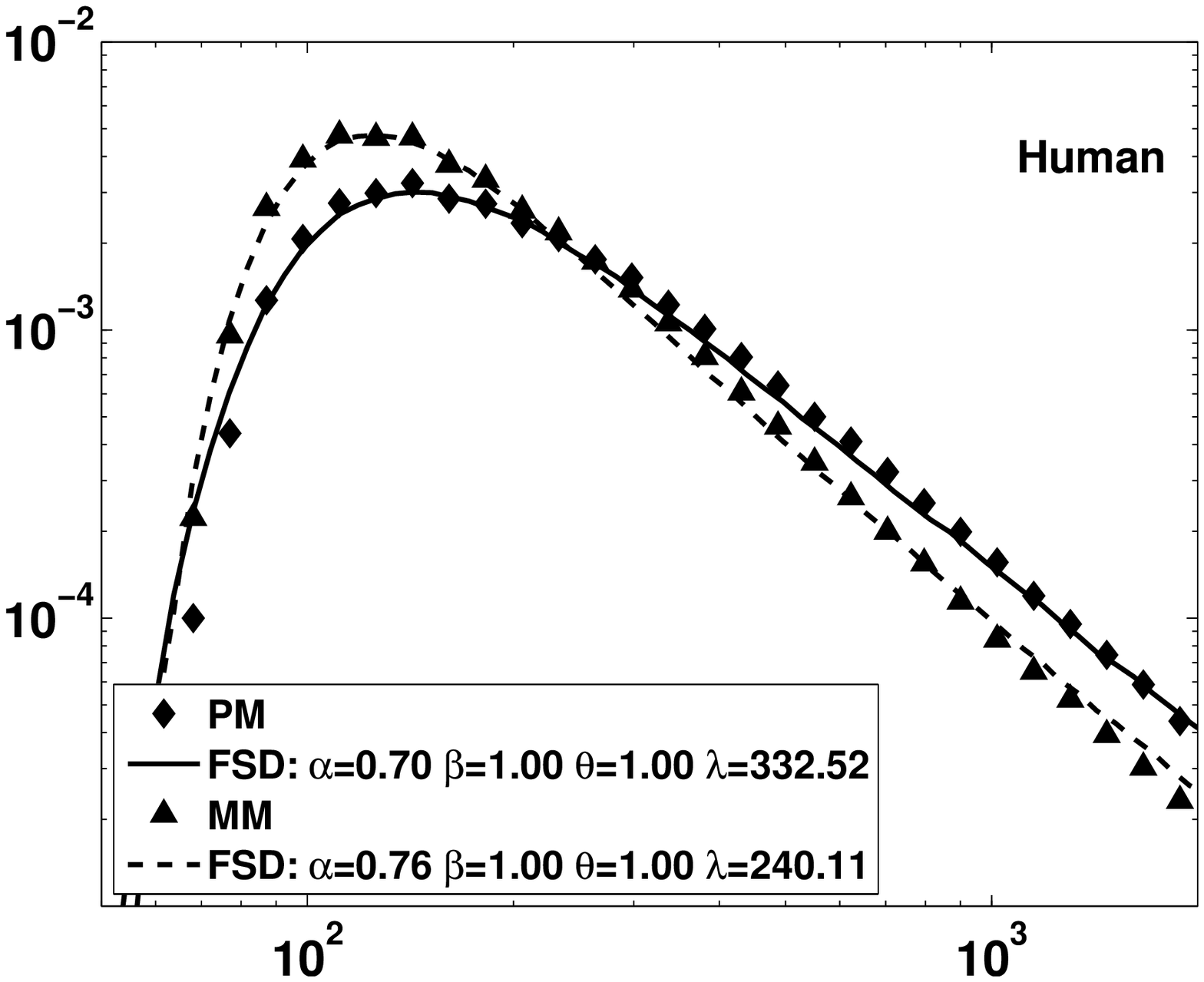}\hfill
\includegraphics[width=0.24\textwidth]{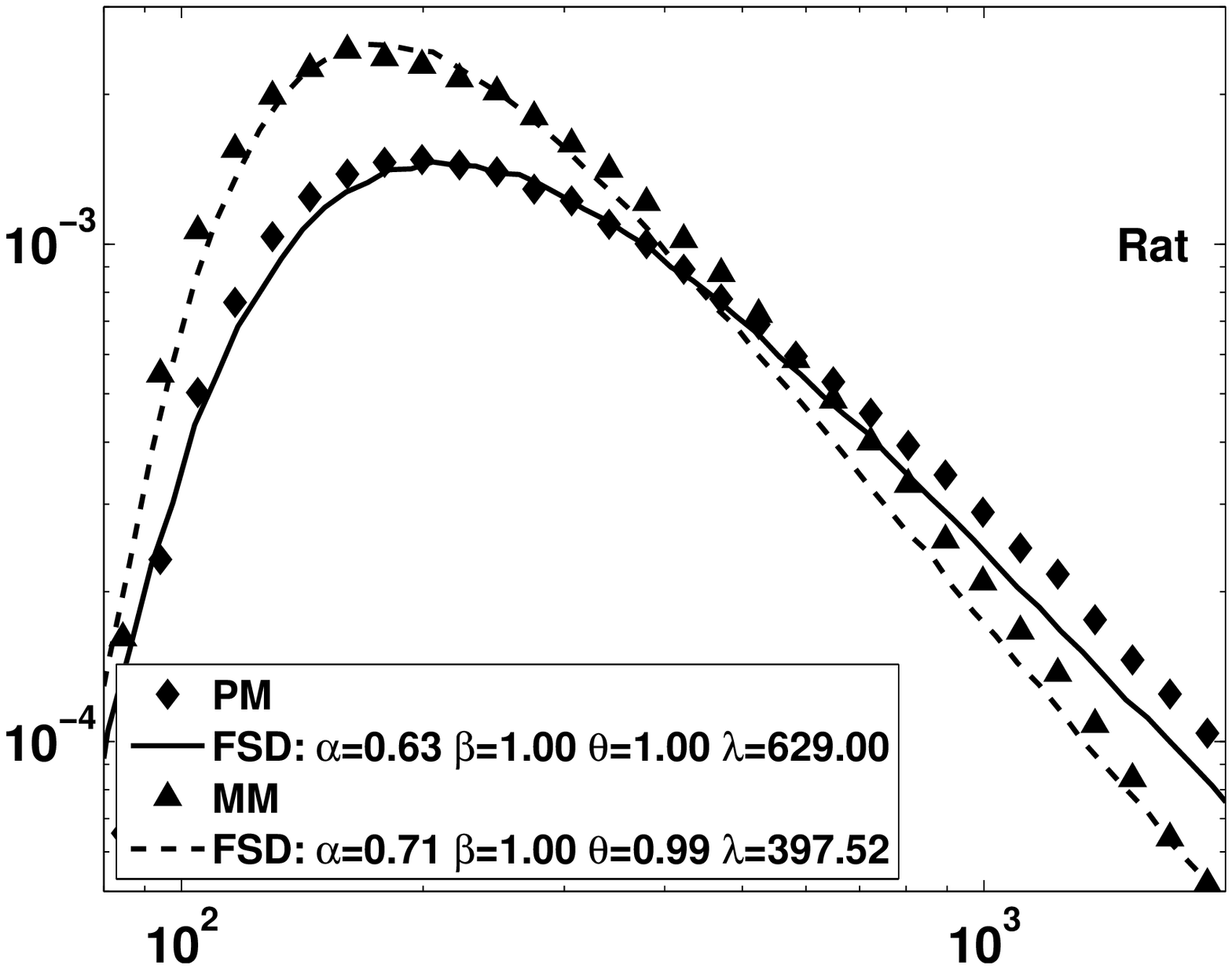}\\
\includegraphics[width=0.24\textwidth]{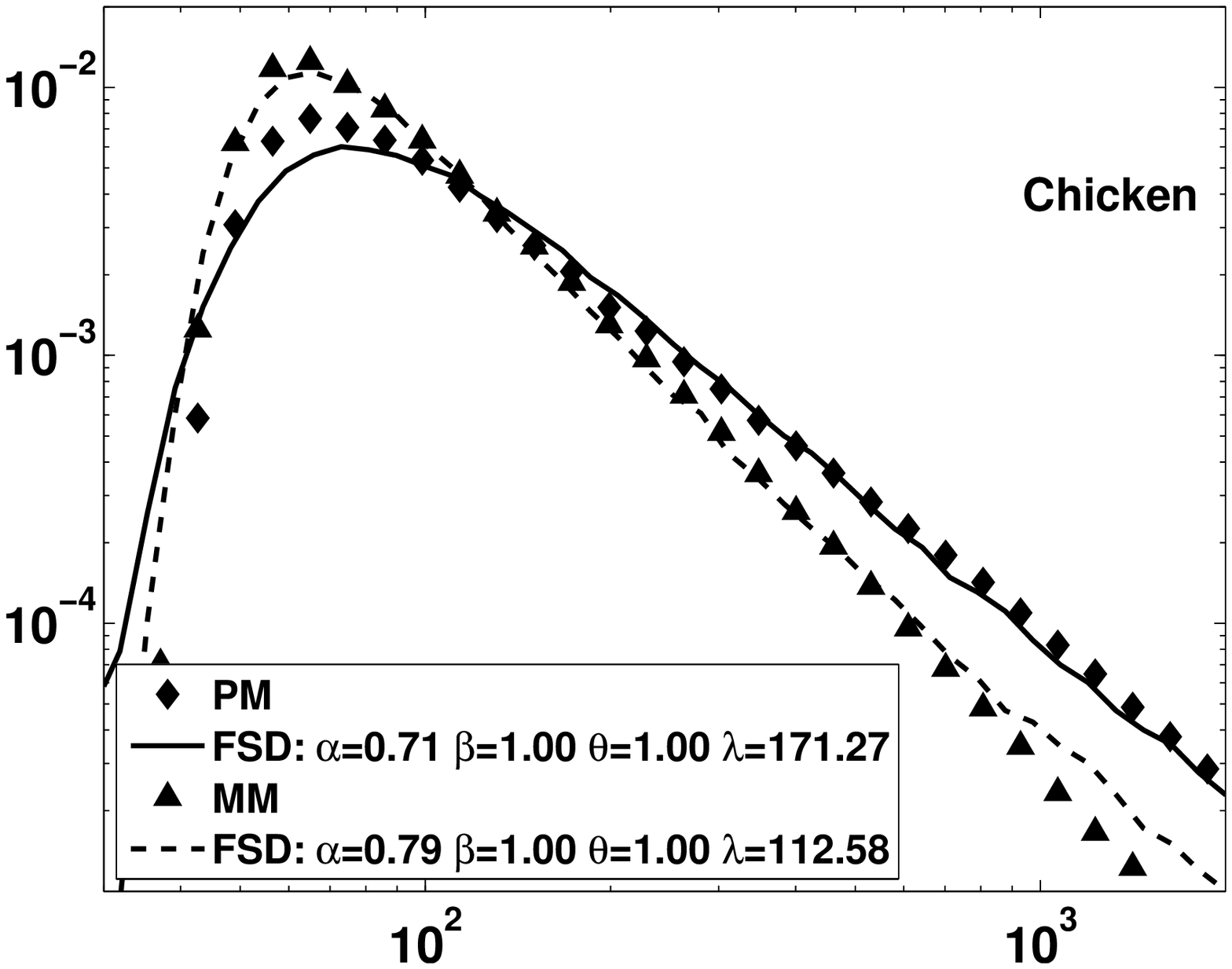}\hfill
\includegraphics[width=0.24\textwidth]{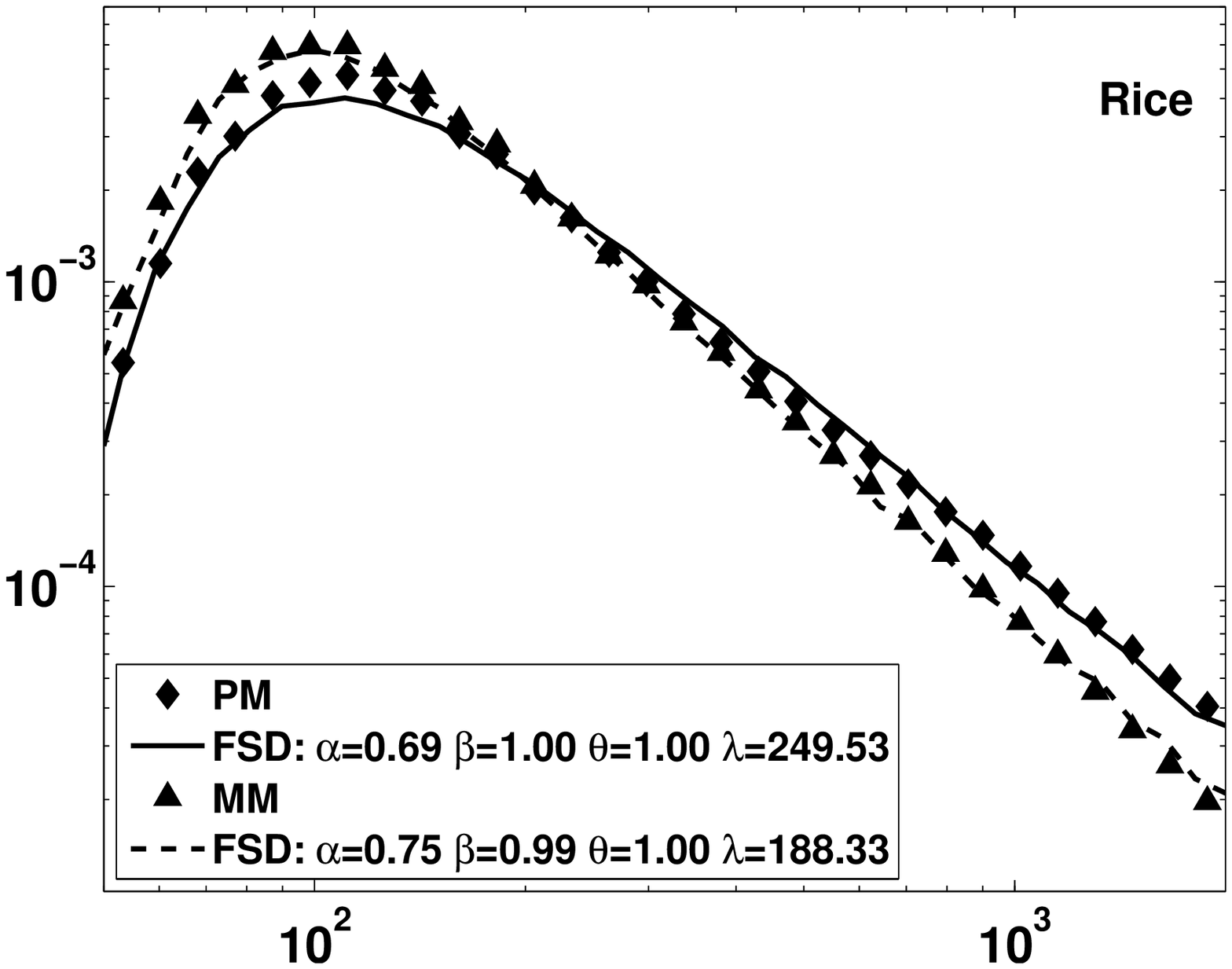}\\
\includegraphics[width=0.24\textwidth]{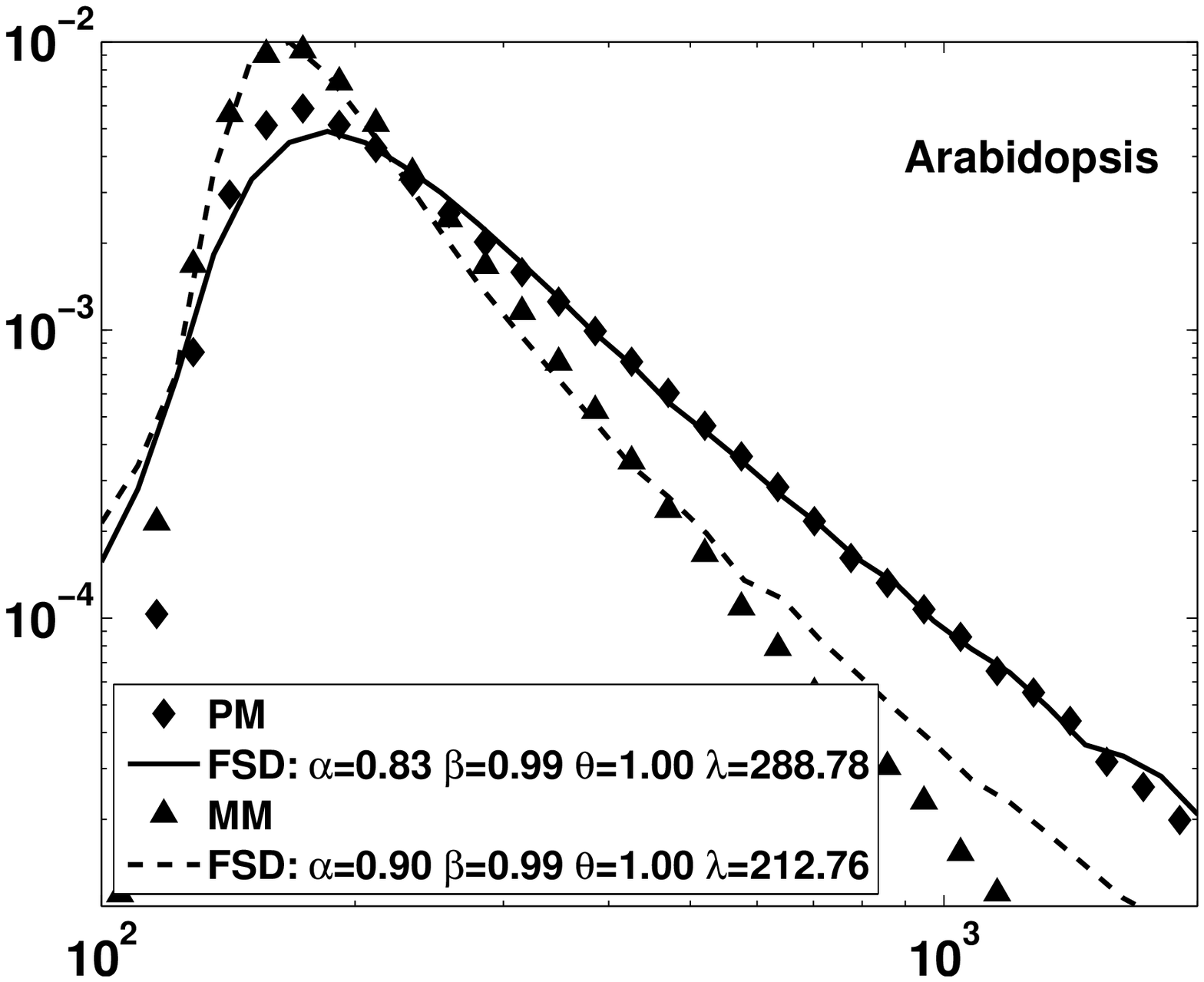}\hfill
\includegraphics[width=0.24\textwidth]{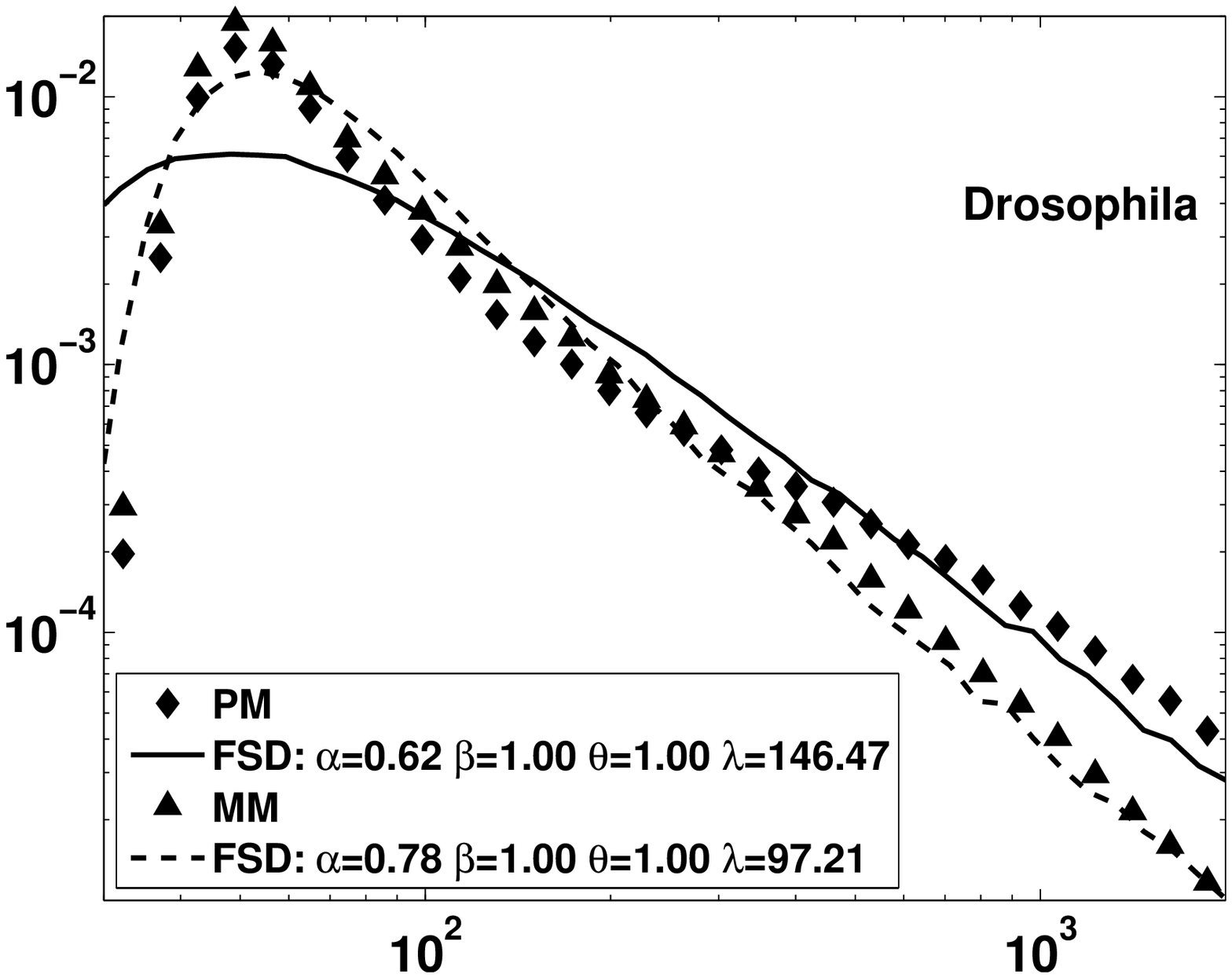}\\
\includegraphics[width=0.24\textwidth]{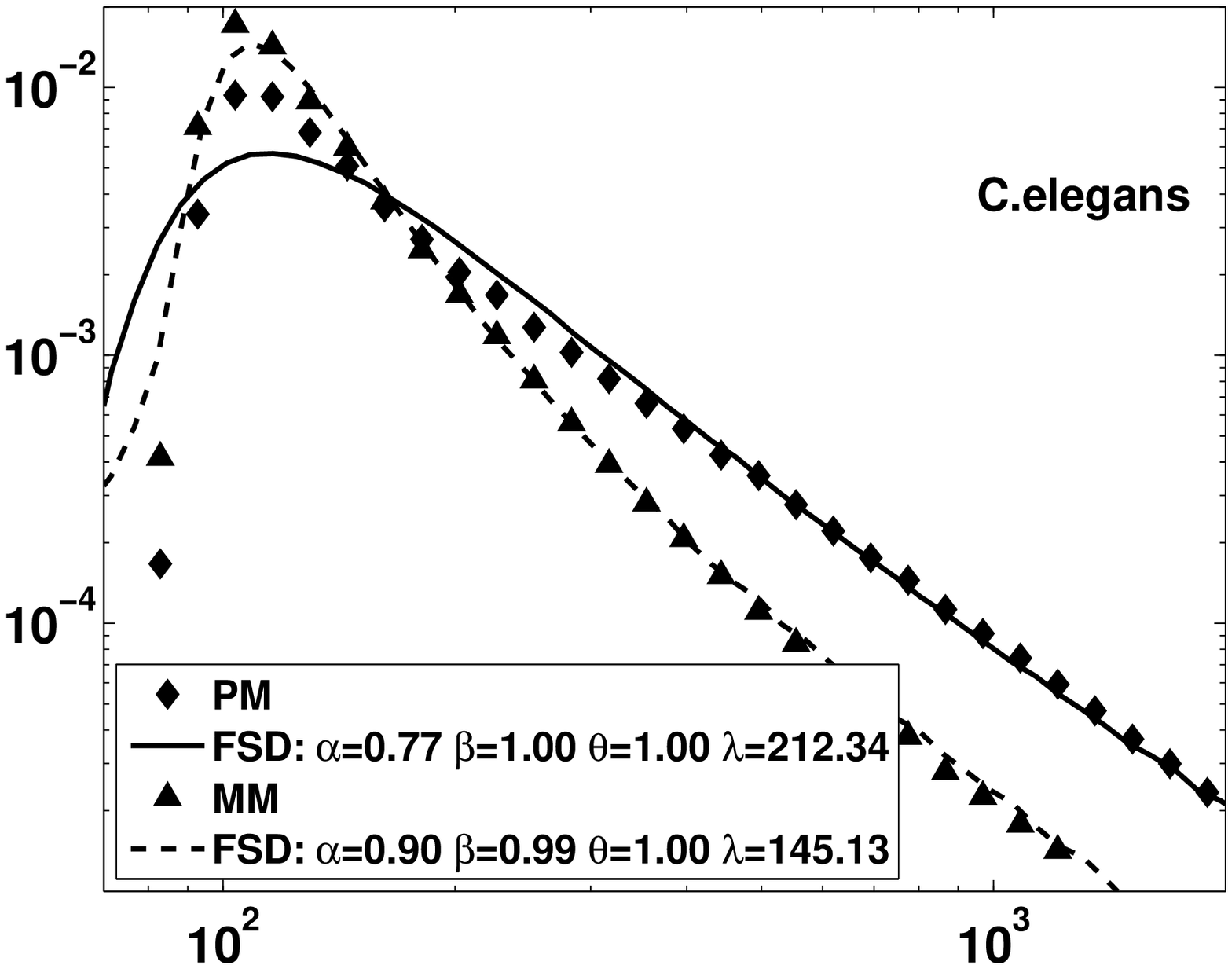}\hfill
\includegraphics[width=0.24\textwidth]{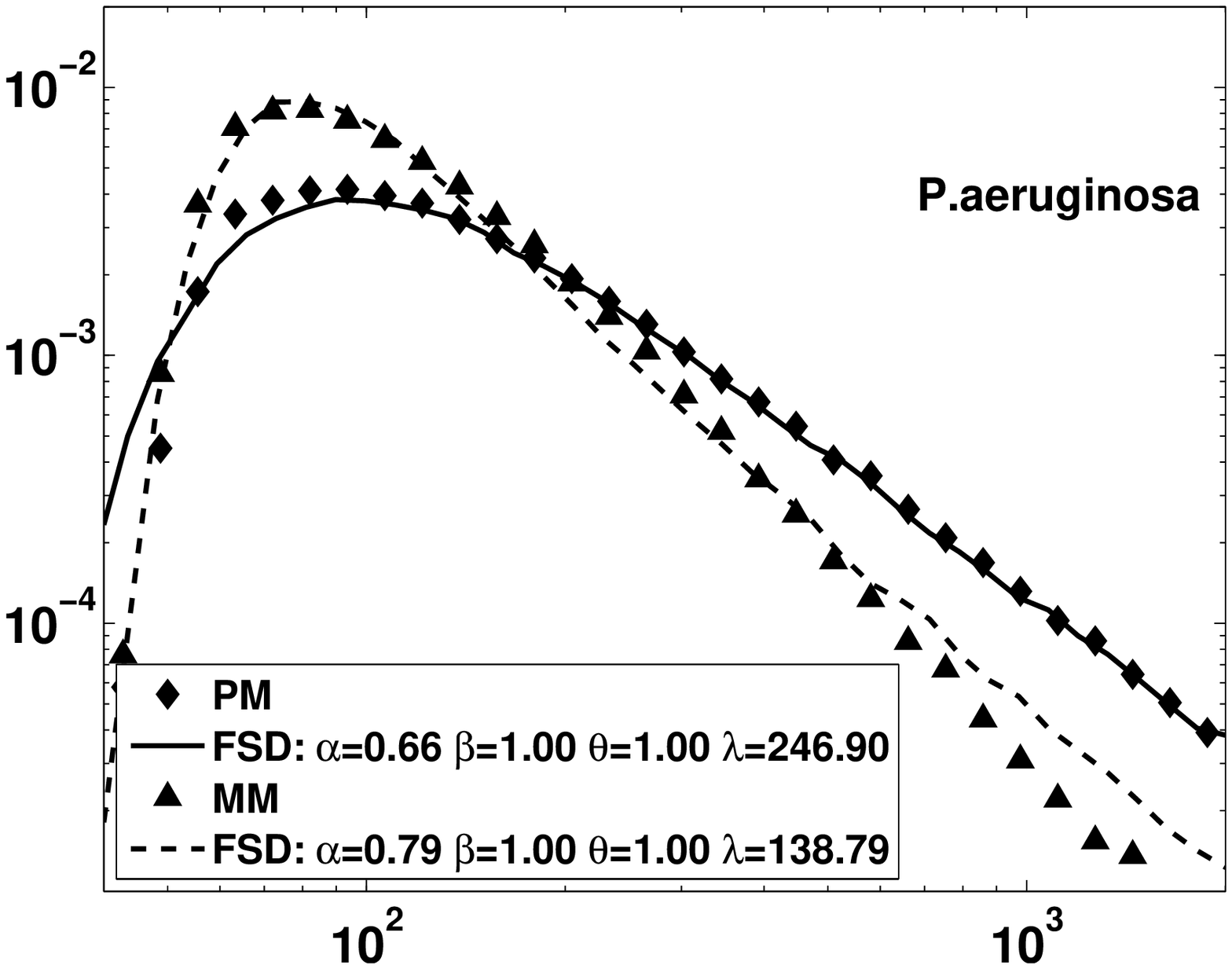}\\
\caption{Distribution of gene expression of microarray experiments for various organisms. Names of organisms are showed on the pictures. Black points are empirical distribution and solid curve is fractional stable distribution. Parameters of FSD are showed on the figures.}\label{fig:Affy_fsdapprox}
\end{figure}

The results of approximation of gene expression profiles for microarrays of the Affymetrix company are shown on the fig.~\ref{fig:Affy_fsdapprox}. On the figure the probability density functions (PDFs) for various organisms are depicted. Diamonds and crosses correspond to PDFs of gene expression for PM and MM channels respectively. Solid line and dashed line correspond to the FSDs are calculated for estimated values of the parameters for PM and MM channels respectively. It is seen from the figure  more satisfactory agreement is achieved for a gene expression of a human, a rat,a chicken and a rice both for PM channel and for MM channel. For C. Elegans and P. aeruginosa an satisfactory agreement between theoretical and empirical distributions is achieved only for MM channel. However when testing the hypothesis of acceptance of two distributions the $\chi^2$ criterion rejects the hypothesis about fractional stable nature distribution of gene expression for all processed organisms. For others results which depicted on fig.~\ref{fig:Affy_fsdapprox} difference of empirical and theoretical distributions are clearly seen.

\begin{figure}
\includegraphics[width=0.23\textwidth]{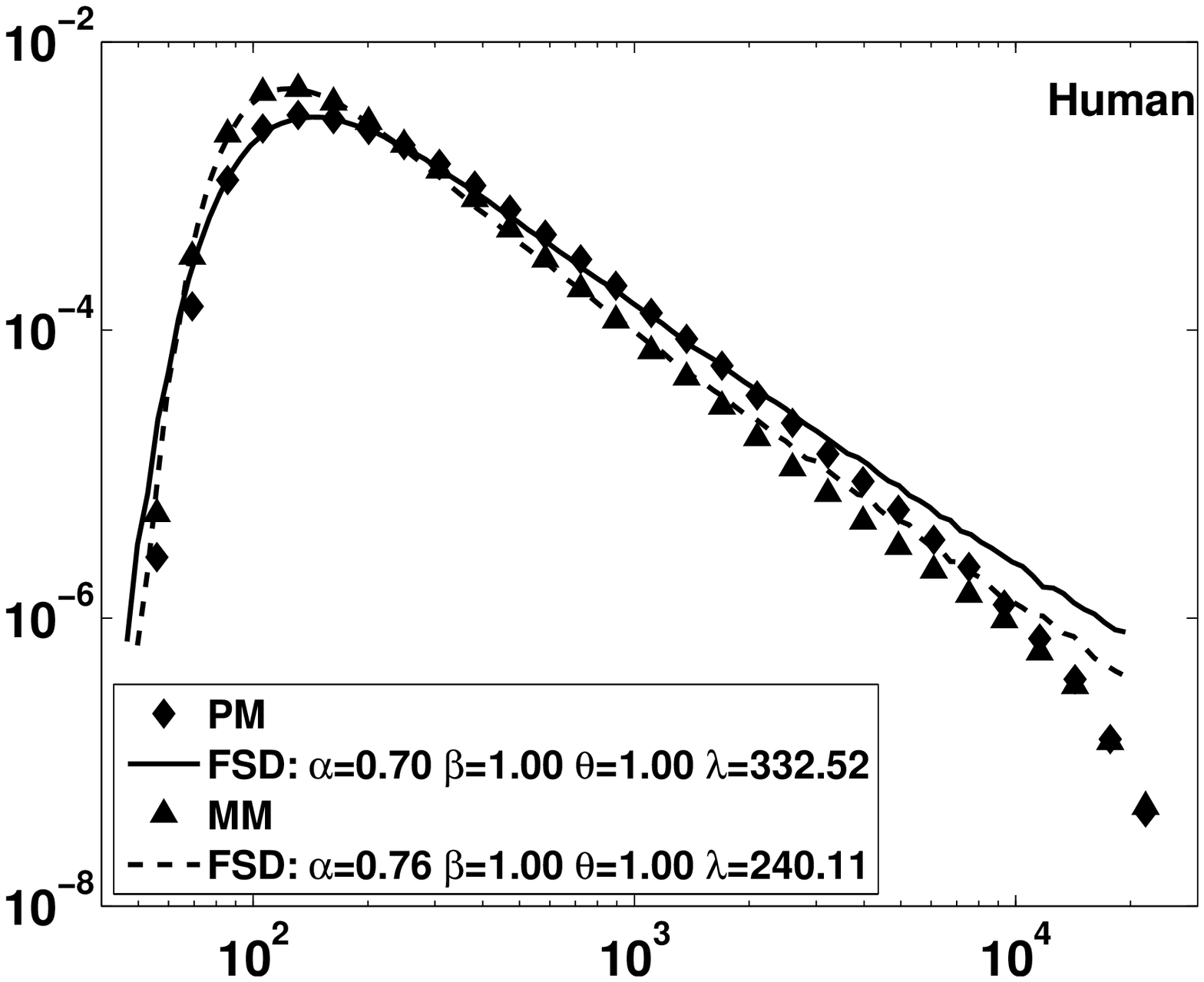}\hfill
\includegraphics[width=0.23\textwidth]{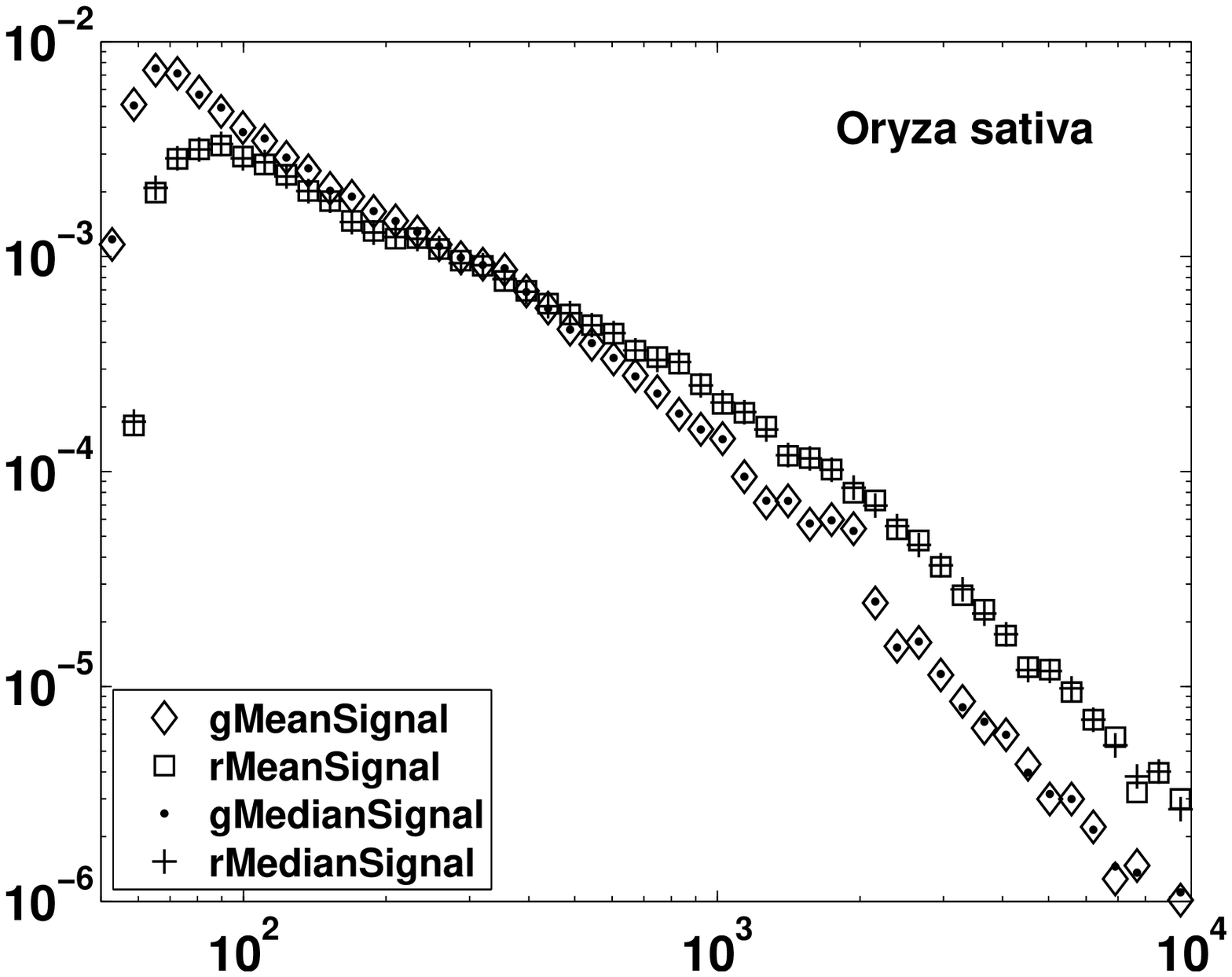}\\
\parbox[t]{0.23\textwidth}{\caption{PDF of gene expression for a human genome obtained by microarray of the Affymetrix company. Notation same as on the fig.~\ref{fig:Affy_fsdapprox}}\label{fig:human_nolim}}\hfill
\parbox[t]{0.23\textwidth}{\caption{PDF of gene expression obtained by microarray of Agilent company. Diamonds (points) are mean (median) signal of green channel, squires (crosses) are mean (median) signal from red channel.}\label{fig:agilen_4ch}}
\end{figure}

Nevertheless it should be noted what this difference may be consequences both  of hardware restriction and imperfection of algorithms selection of point glow and their digitization during process of translating them from image to a data file.  One evidence of the presence of hardware constraints may serve fig.~\ref{fig:human_nolim}. On this figure gene expression of human genome is depicted but at the same the empirical distribution has been plotted in all range of values. Here it should be noted what on the fig.~\ref{fig:Affy_fsdapprox} PDFs are plotted not  for all range of gene expression. It is seen from the fig.~\ref{fig:human_nolim} at large values of expression $\gtrsim 10^4$ a power law dependence is broken and PDF rapidly goes to zero. Such effect is called an effect of truncation and may be consequence of the hardware restriction at large values of gene expression intensity.

Let consider now the results of processing microarrays of the Agilent company.  In the RAW files  four channels correspond to gene expressions results. These channels differ by color and by the method of calculation of gene expression. In technological process of these microarrays red and green dye are used and two method of calculation of gene expression value are also used. The first method consists in calculation of mean value of intensity obtained from all pixels a probe under investigation. The second method consists in choosing median value of intensity of gene expression at processing of all pixels of the probe. According to this here and after we will denote: gMeanSignal (rMeanSignal) is mean signa of green (red) channel; gMedianSignal (rMedianSignal) is median signal in green (red) channel. During the process of processing it was obtained what PDFs of median and mean signal from same color almost coincide with each other.  It is clearly seen from the fig.~\ref{fig:agilen_4ch} on which PDFs of gene expression are depicted for a genome of a rice (Oryza sativa). From the figure we can see  that PDFs of gene expression for mean and median signals for both channels practically coincide with each other. Same conclusions were obtained for all processed experimental data. Therefore in this work we will be use only median signal from red and green channels.

\begin{figure}
\includegraphics[width=0.24\textwidth]{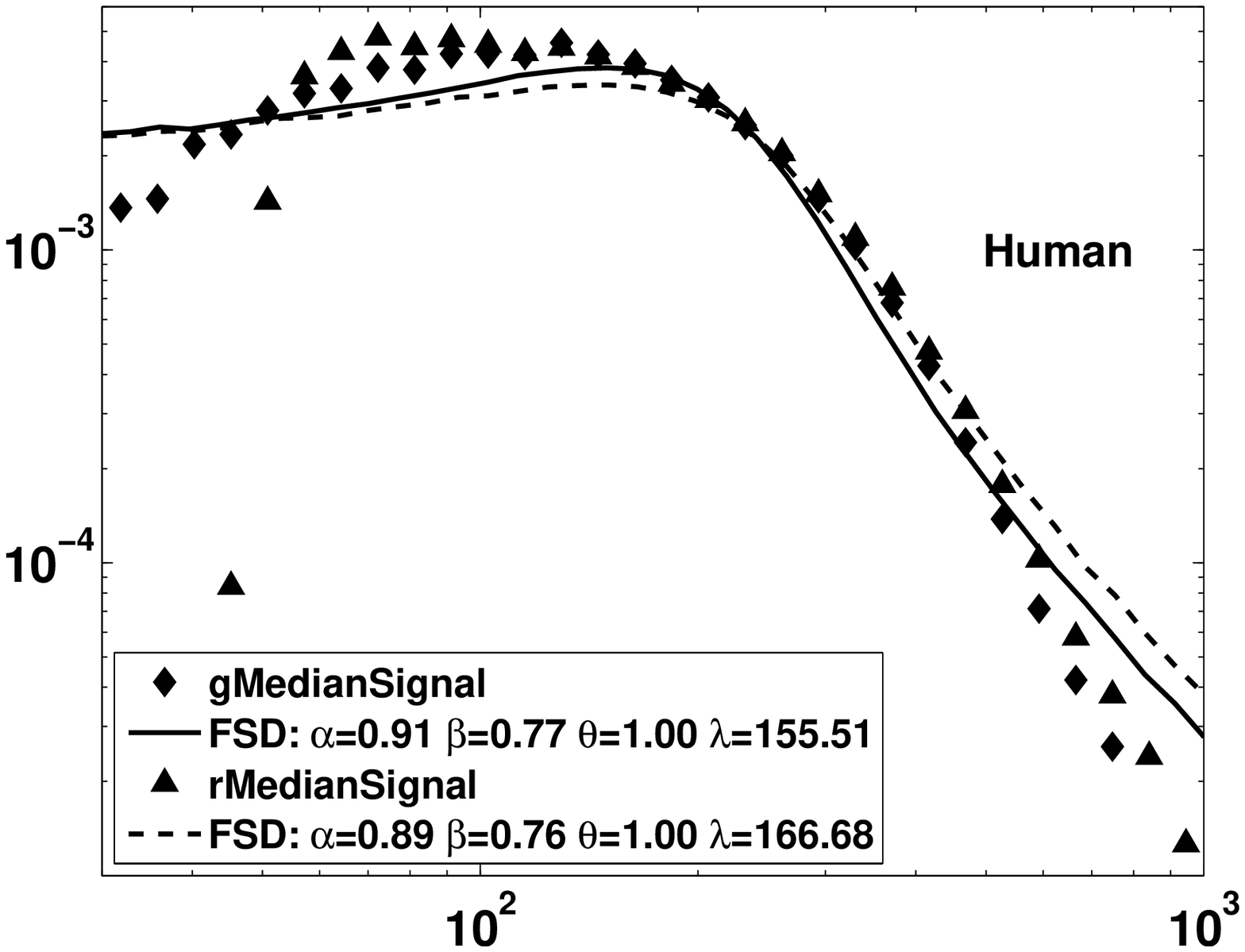}\hfill
\includegraphics[width=0.24\textwidth]{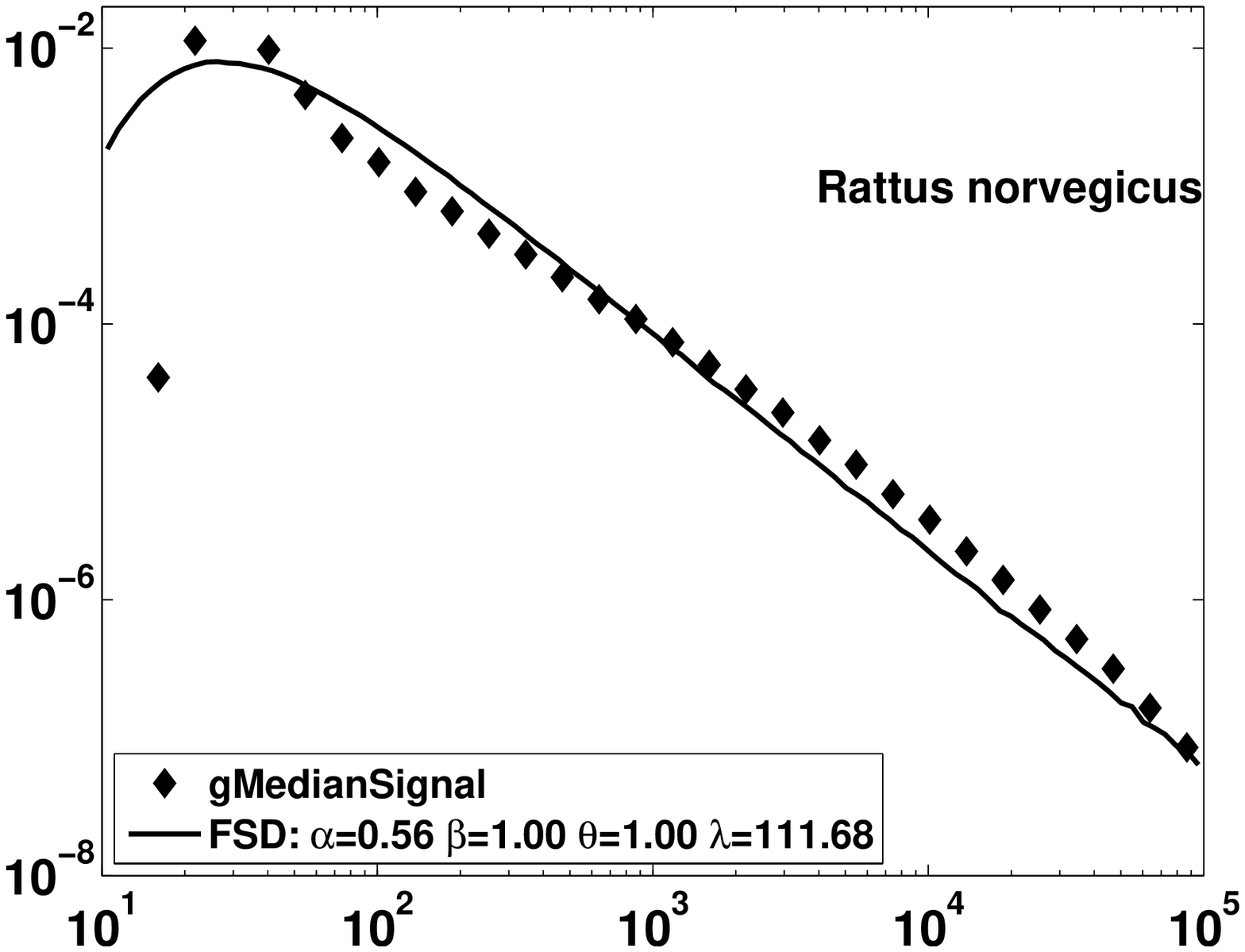}\\
\includegraphics[width=0.24\textwidth]{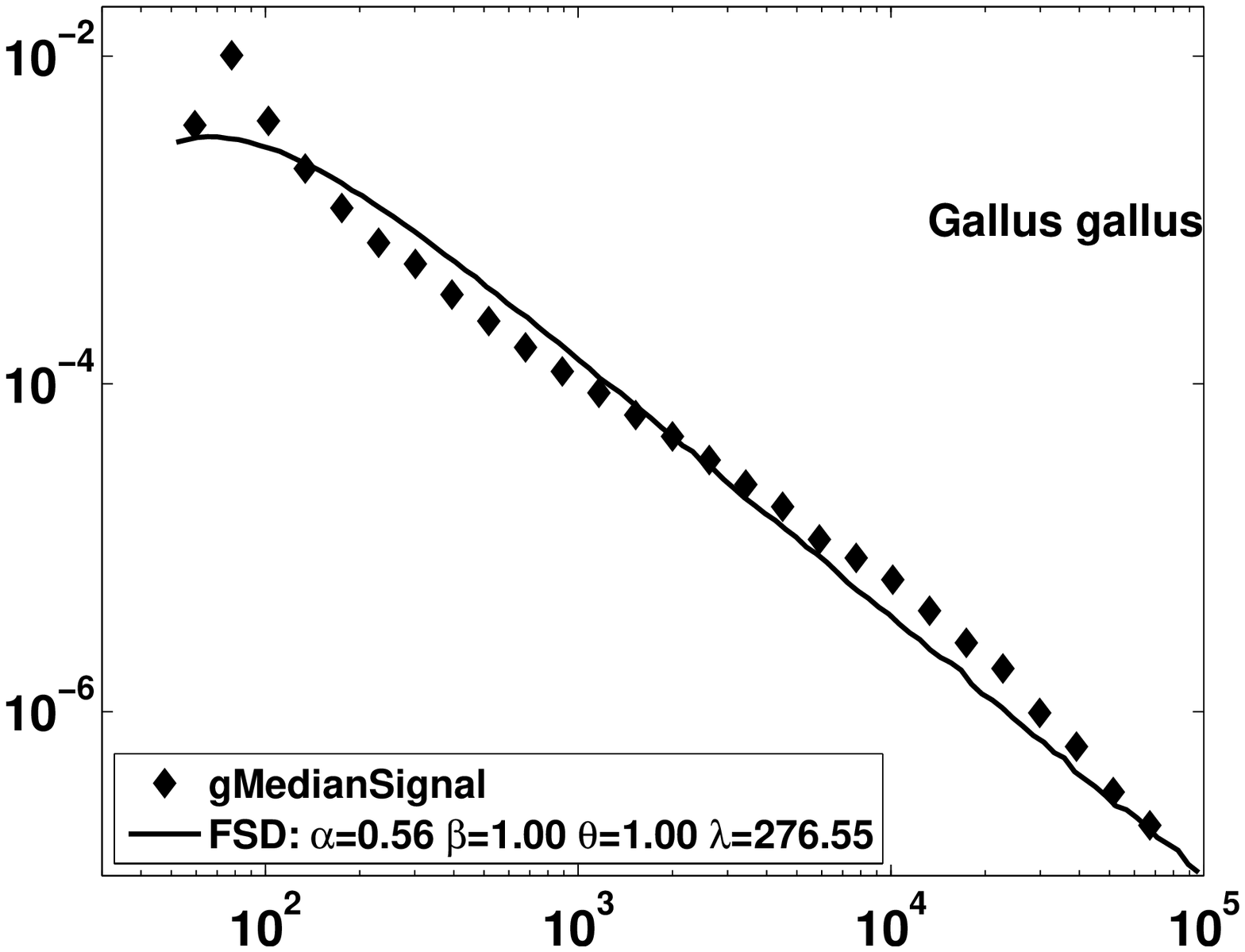}\hfill
\includegraphics[width=0.24\textwidth]{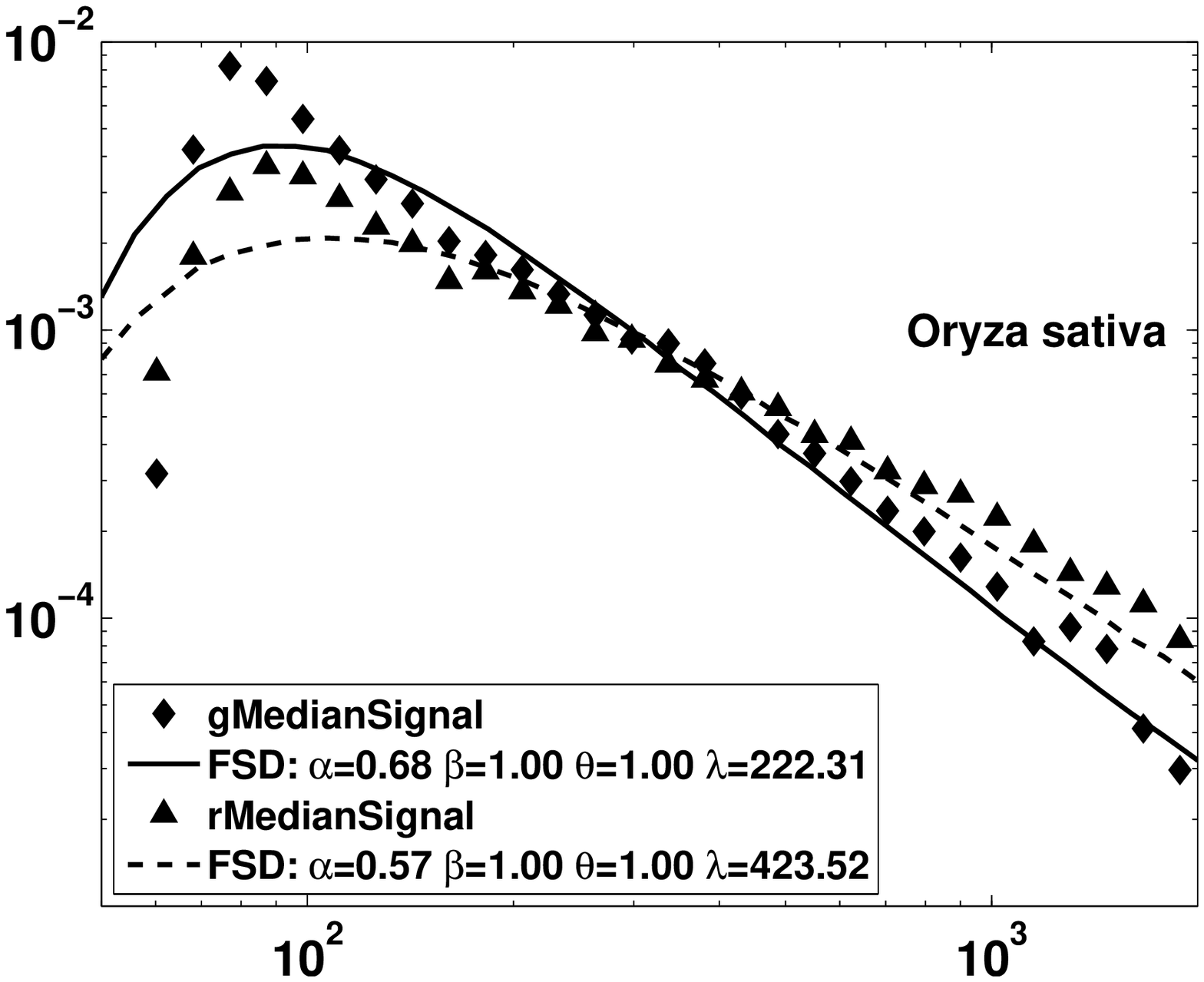}\\
\includegraphics[width=0.24\textwidth]{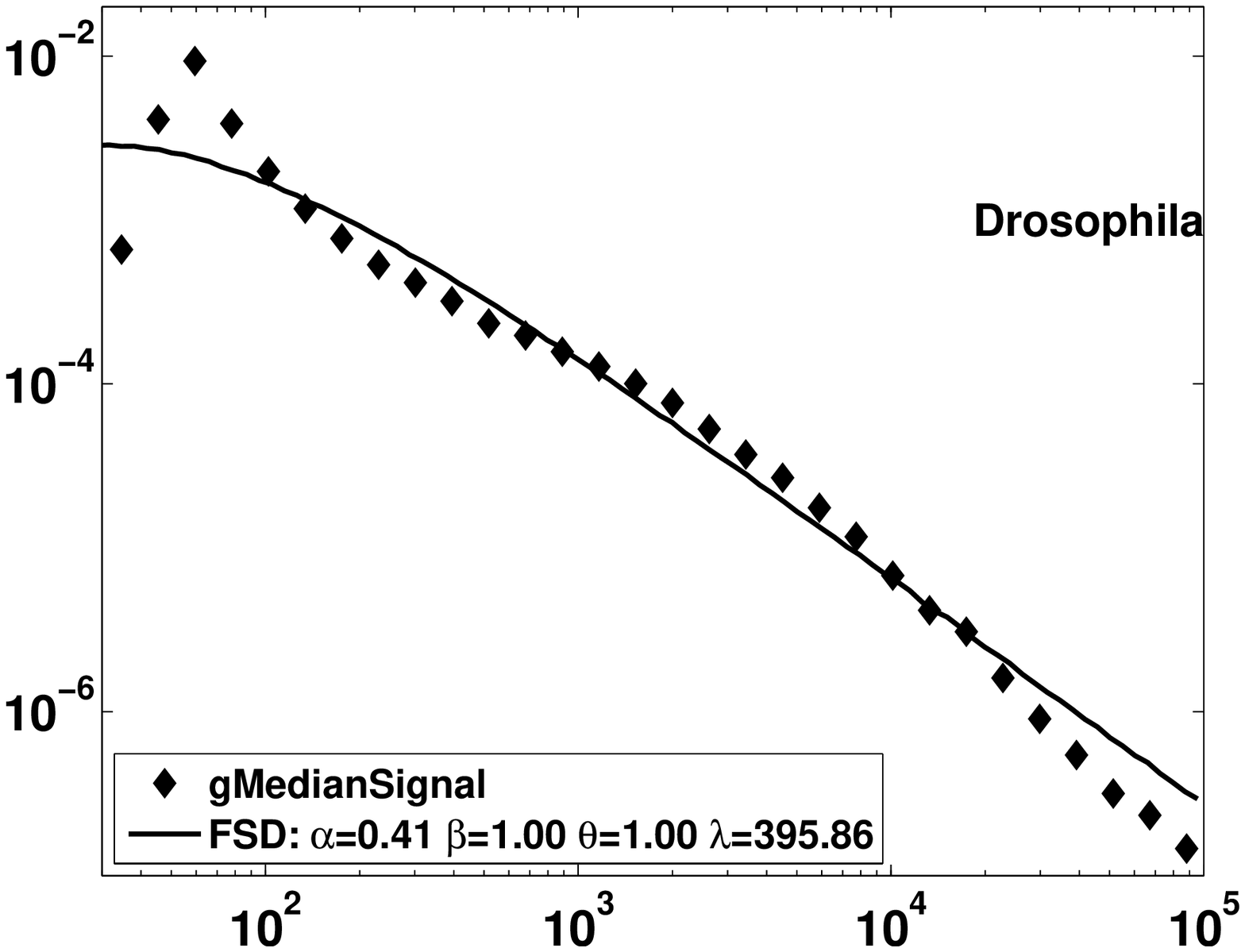}\hfil
\includegraphics[width=0.24\textwidth]{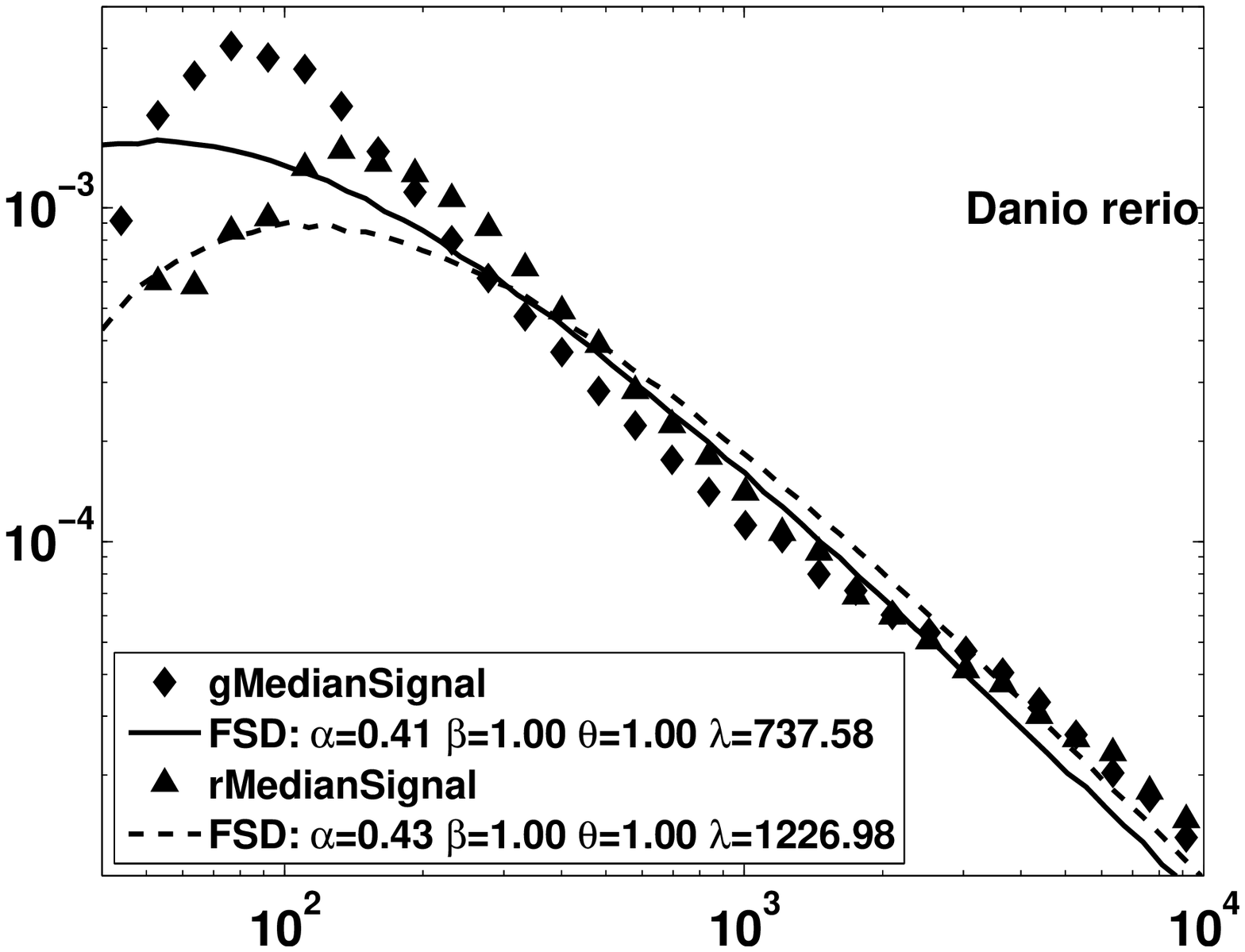}\\
\includegraphics[width=0.24\textwidth]{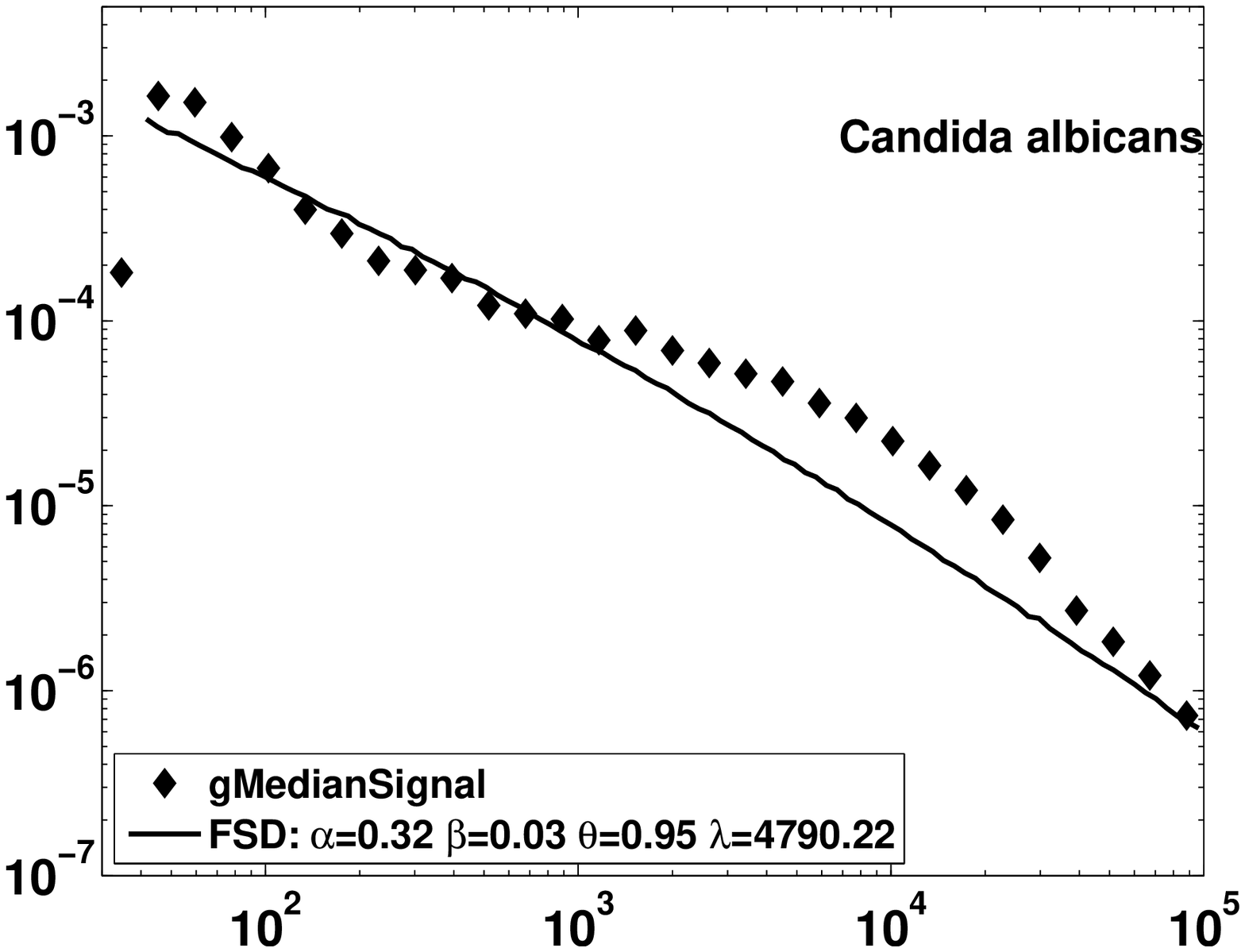}\hfill
\includegraphics[width=0.24\textwidth]{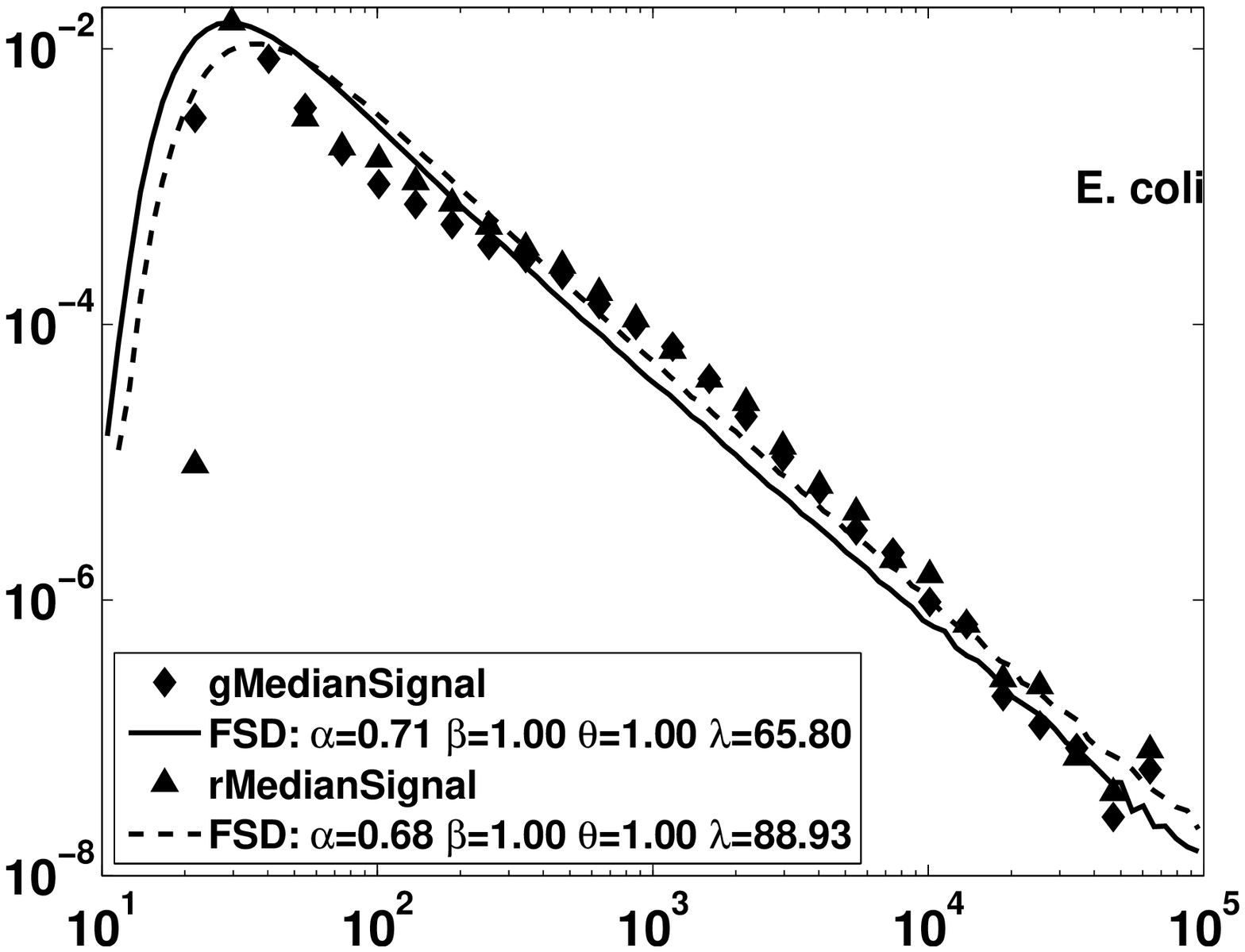}
\caption{Distribution of gene expression of microarray experiments for various organisms obtained from the Agilent microarray chip. Names of organisms are showed on the pictures. Black points are empirical distribution and solid curve is fractional stable distribution. Parameters of FSD are showed on the figures.}\label{fig:Agilent_fsdapprox}
\end{figure}

For microarrays of the Agilent company were selected experimental data for mammals (Homo sapiens, Rattus norvegicus), bird (Gallus gallus), fish (Danio rerio), plant (Oryza sativa), insect (Drosophila melanogaster), fungus (Candida albicans)  and bacterium (E. coli). The empirical PDFs for the median signal from red and green channels and PDFs of FSDs are shown on the fig.~\ref{fig:Agilent_fsdapprox}.  It is clearly seen that empirical PDFs aren't FSDs. Disagreement of empirical and theoretical distributions is very substantially. Nevertheless, let distinguish  some properties  which inherent to all the processed data.  It is clearly seen that the asymptotic of the experimental PDFs haven't power law dependence $\propto x^{-\alpha-1}$. Most likely we can talk about dependence  which close to power-law behaviour. Such behaviour differs from the results obtained by using microarrays the Affymetrix manufacture (see.~fig.~\ref{fig:Affy_fsdapprox}). An existence of hardware distortions and distortions of algorithms of translating of intensity from an image file to numerical value can serve causes of deviation from the power-law dependence.

\begin{figure}
\includegraphics[width=0.24\textwidth]{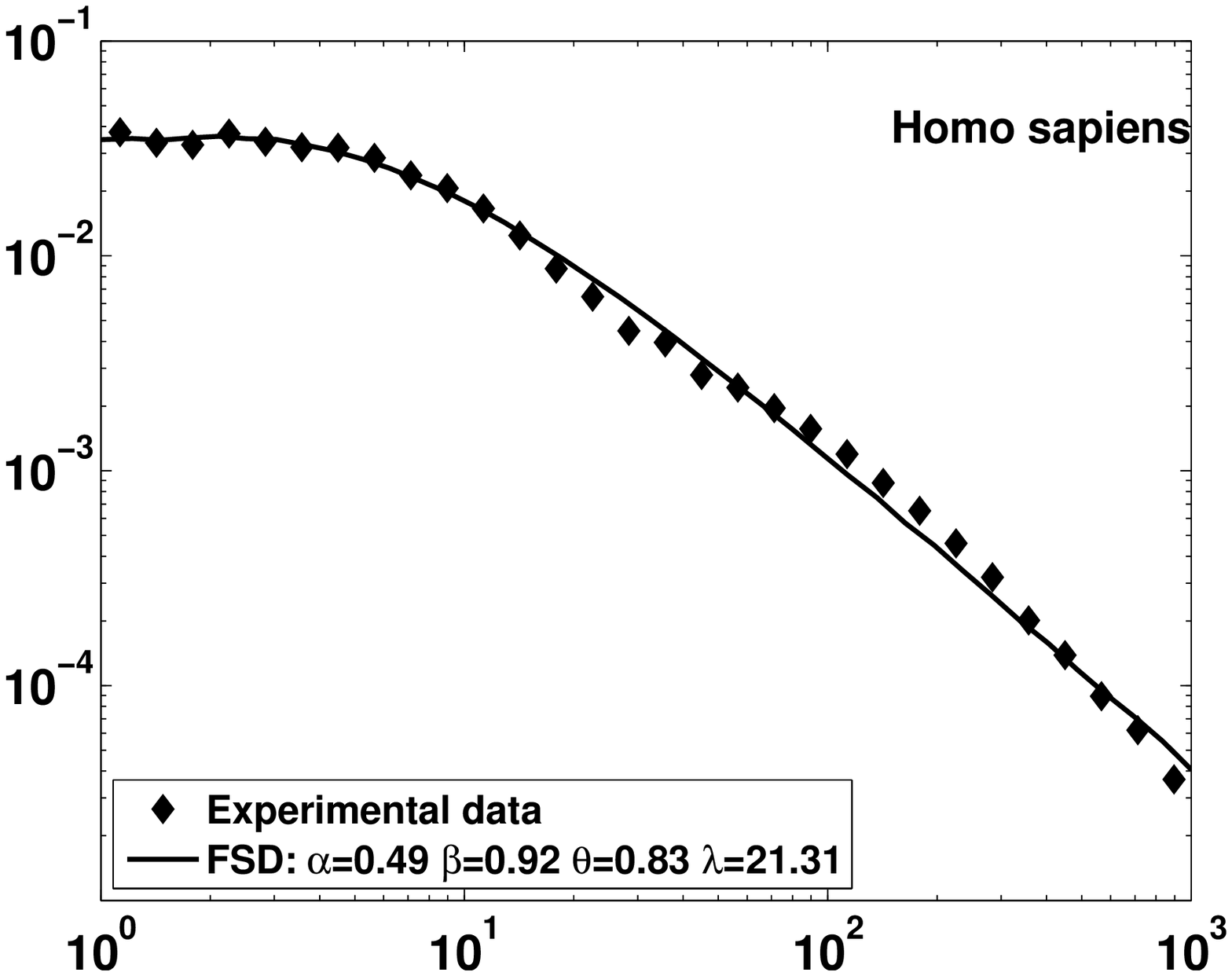}\hfill
\includegraphics[width=0.24\textwidth]{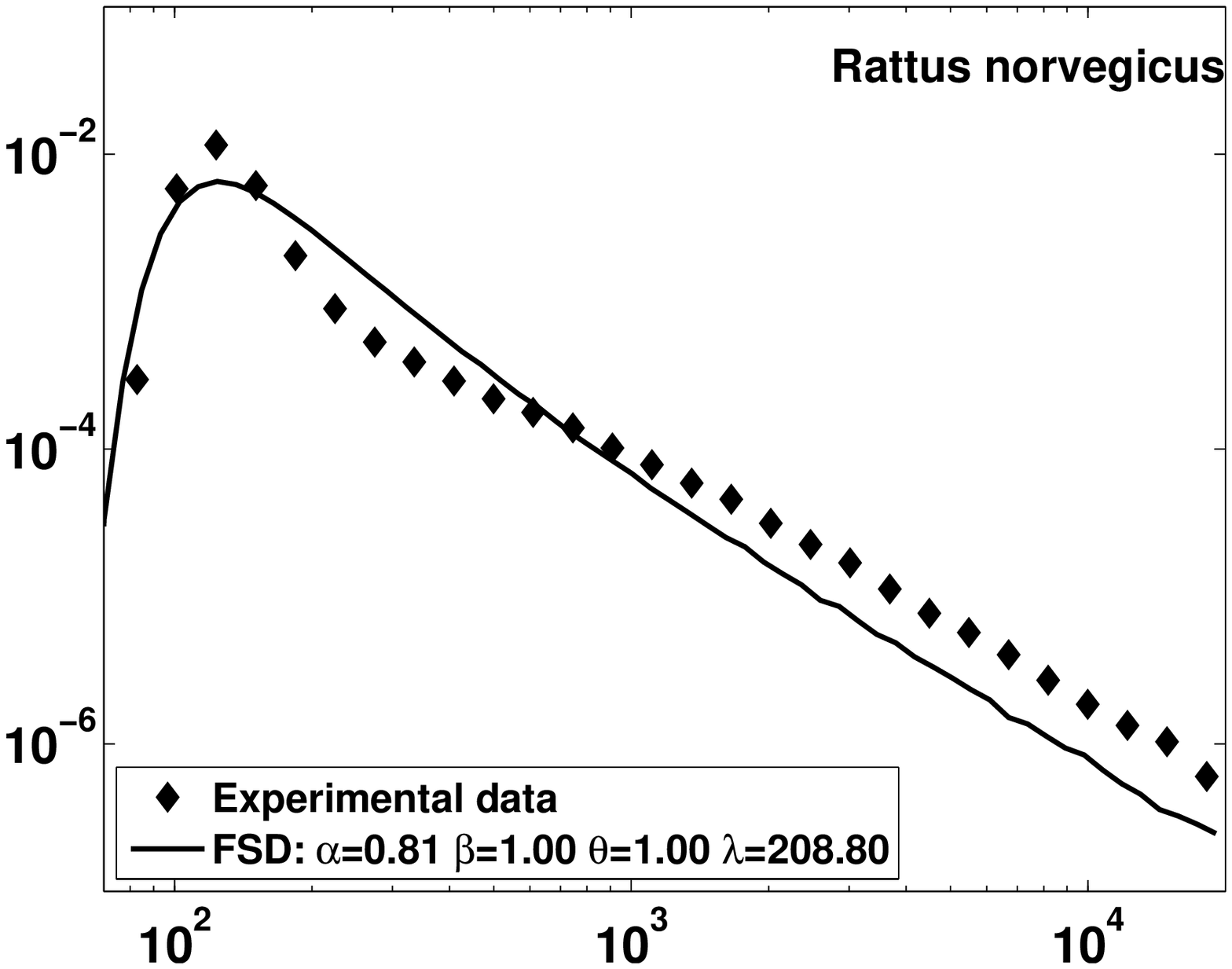}
\caption{Distribution of gene expression of microarray experiments for various organisms obtained from the Illumina microarray chip. Names of organisms and the estimated parameters on the FSDs are showed on the pictures. Diamonds are empirical distribution and solid curve is fractional stable distribution. Parameters of FSD are showed on the figures. }\label{fig:Illumina_fsdapprox}
\end{figure}

The PDFs of gene expression of human (Homo sapience) and rat (rattus norvegicus) for Illumina microarrays are shown on the fig.~\ref{fig:Illumina_fsdapprox}. On the figures diamonds are experimental PDFs and solid line are FSD. It is seen from the figures the satisfactory agreement between experimental and theoretical PDFs is observed only for human genome. However usage $\chi^2$ Pirson's criterion for checking correctness   the hypothesis about fractional stable nature of the experimental distributions leads to necessity to reject this assumption. Nevertheless, it is seen from the figure that the FSD is good approximation for PDF of gene expression for human genome. For another genome is presented here the experimental distribution aren't belong to the class of FSDs. As well as in the previous case the asymptotic of the experimental PDFs isn't described by power-law dependence. The power-law dependence is observed in mean but on this dependence some distortions are imposed.

\section{Results and Discussion}

In present work the attempt was made to approximate distribution of gene expression by FSD. It is necessary to know four parameters for unique determination of the FSD. Therefore the one of the main tasks which has been solved here is the task of estimation of the parameters of the FSD by sample of independent identical distributed random variable. The estimation algorithm is described in the Section~\ref{sec:FSDEst}. Next by estimated values of the parameters the histograms of FSD was being constructed and these distribution were compared with experimental histogram. The $\chi^2$ Pearson's criterion was applied for test the hypothesis about coincidence of two distributions.

An object of investigation  were selected several organisms belonged to various classes: mammals, birds, fishes, plants,  fungus, bacterium. Since in the present time there are many  companies which produce microarrays, then the interesting question appears: how relate PDFs of gene expression between each other which have been obtained by microarrays of various companies? From fundamental understanding it is clear; since genes expression is proportional their concentration then law of distribution must be invariant towards manufacturer of microarray platform. In this work microarrays of three manufacturers (Affymetrix, Agilent è Illumina) were selected.

The results of comparison of theoretical end empirical densities are presented on the figs.~\ref{fig:Affy_fsdapprox},~\ref{fig:Agilent_fsdapprox},~\ref{fig:Illumina_fsdapprox}. As seen from the figures the law of distribution of gene expression for microarrays of different manufacturers is different. The PDFs of gene expression for microarrays the Affymetrix manufacture have clearly marked power-law asymptotic. However the effect of truncation is observed at large value of intensity of gene expression. (see~fig.~\ref{fig:human_nolim}) which breaks the power-law asymptotic. Clearly marked power-law asymptotic doesn't observes for PDFs of gene expression for microarrays the Agilent and Illumina manufacture (see~figs.~\ref{fig:Agilent_fsdapprox}~è~\ref{fig:Illumina_fsdapprox}). Is observed some decreasing which resembles  the power-law dependence. Therefore for these data we can't talk about power-law asymptotic. By our opinion the differences in used algorithms of processing of initial data at their reading from microarray and subsequent translating there from image file to a numerical value are causes of divergence   between the results of different platforms.

Approximation of PDFs of gene expression by FSDs has showed that the best agreement is achieved for gene expression of mammals and plants for microarrays the Affymetrix manufacture (see~fig.~\ref{fig:Affy_fsdapprox}). However $\chi^2$ criterion rejects hypotheses about coincidence  of these two distributions. For PDFs of gene expression of microarrays the Agilent and Illumina manufacture the situation is absolutely different. There is clear difference between experimental distributions here and FSDs and in this case we can't talk about of coincidence of these distributions.

Nevertheless, the FSD good enough approximates empirical distribution both in the central part and in the tail part for gene expression of mammals and plants genomes. As we can see the values of the parameter  $\alpha$ lie within interval $0.62\leqslant\alpha\leqslant0.83$. This values are in good agreement with results of works \cite{Ueda2004a,Furusawa2003,Kuznetsov2002}. A value of second characteristic parameter of FSD $\beta$ little differs from unit. This means that distribution of gene expression belongs to the class of stable laws which is a subclass of FSDs. As we can see the FSD good approximate empirical data of gene expression.

The fact that empirical distribution of gene expression is described by FSD allows to make some assumption about character of background processes. As was noted above, the FSD is the limit distributions of sums of independent identically distributed random variables. Physical interpretation of the sum (\ref{eq:sumFSD}) is a trajectory of particle undergoing a random walk. In this process, random variables $X_{ij}$ is random races  and $T_{ij}$ have mean random rest time between two successive jumps in $i$-th trajectory. Thus, sums (\ref{eq:sumFSD}) and (\ref{eq:sumTime}) describe process of random walks with instantaneous jumps. Such process named Continuous Time Random Walk (CTRW) \cite{Metzler2000}. In the work \cite{Uchaikin2000} was shown that limit distribution of particle coordinate in framework of CTRW process is expressed through FSD. As consequence we can assume  that the basis of the processes leading to the observed distribution of gene expression levels, are the processes described scheme CTRW.

On the other hand it is known that asymptotic behavior of CTRW process is described by generalized diffusion equation \cite{Metzler2000} expressed through fractional derivatives
\begin{equation}\label{eq:fracdiffeq}
\frac{\partial^\beta p(x,t)}{\partial t^\beta}=-D\left(-\Delta\right)^{-\alpha/2}p(x,t)+\frac{t^{-\beta}}{\Gamma(1-\beta)}\delta(x).
\end{equation}
Here $\partial^\beta/\partial t^\beta$ is Riemann-Liuville fractional derivative and  $(-\Delta)^{-\alpha/2}$ is Laplace operator of fractional order \cite{Samko1973}, $D$ is diffusion constant. Solution of this equation is expressed through FSD \cite{Uchaikin2000}
$$
p(x,t)=\left(Dt^\beta\right)^{-1/\alpha}q\left(|x|\left(Dt^\beta\right)^{-1/\alpha};\alpha,\beta,0,1\right),
$$
where $q(x;\alpha,\beta,\theta,\lambda)$ is FSD (\ref{eq:FSD}). At the same time the parameters $\alpha$ and $\beta$ simultaneously are exponents of fractional power of derivatives in the equation (\ref{eq:fracdiffeq}). Thus, from this facts, we can conclude, that processes leading to observed gene expression can be described by using equation in fractional derivatives. But the question about nature and main characteristics of these processes remains open.

\section*{Acknowledgement}

\paragraph{Funding\textcolon} The work was supported by Ministry of Education and Science of Russian Federation (grant No 2014/296, No 6.1617.2014/K) and Russian Foundation or Basic Research (grant No. 12-01-00660)

\appendix
\section{Simulation of Fractional Stable Random Variables}\label{app:FSRV}

According to the work \cite{Kolokoltsov2001} FS random variable can be can be represented as ratio of two strictly stable random variable (\ref{eq:fsrv}).
For simulating $Y(\alpha,\theta)$ the Chamber's algorithm~\cite{Chambers1976}
\begin{align*}
Y(\alpha,\theta)&= \lambda^{1/\alpha}\sin(\alpha(V+C_1))(\cos V)^{-1/\alpha}\times\\
&\times\left(\cos(V-\alpha(V+C_1))/W\right)^{(1-\alpha)/\alpha},\quad \alpha\neq1\\
Y(1,\theta)&=(\pi/2)\lambda\tan V,\quad \alpha=1.
\end{align*}
was used, where $C_1=\alpha\theta/(\alpha-1-\mathrm{sign}(\alpha-1))$, $V=\pi(0.5-U_1)$, $W=-\log U_2$. The random variable $S(\beta,1)$ simulated according to Kanters's algorithm \cite{Kanter1975}
$$
S(\alpha)\stackrel{d}{=}
\frac{\sin(\alpha\pi U_3)[ \sin((1-\alpha)\pi
U_3)]^{1/\alpha-1}}{[\sin (\pi U_3)]^{1/\alpha}[-\log U_4]^{1/\alpha-1}},
$$
where $U_1$, $U_2$, $U_3$ and $U_4$ are variables uniformly distributed in $(0,1]$.

\section{Estimation of the parameters by moment method} \label{sec:prmestmom}

Let $Z_1,Z_2,\dots,Z_n$, $n\leqslant4$ be independent, identically distributed random variables with density (\ref{eq:FSD}). The problem is to determine estimates $\hat\alpha, \hat\beta, \hat\theta, \hat\lambda$ of unknown parameters $\alpha,\beta,\theta,\lambda$. This problem was solved in \cite{Bening2004}, where a factional stable stochastic variable was represented in the form (\ref{eq:fsrv}).

Here, we only present the final result. The formulas for estimates $\hat\alpha$, $\hat\beta$, $\hat\theta$, $\hat\lambda$ of parameters $\alpha,\beta,\theta,\lambda$ has the form
$
\hat\theta=1-\frac{2}{n}\sum_{j=1}^n\mathbf{I}(Z_j<0),\
\hat\alpha=\frac{2\pi}{\sqrt{12V_n+\pi^2\left(2Z_n+3\hat\theta^2-1\right)}},\
\hat\beta=A_n\hat\alpha,\
\hat\lambda=\exp\left\{U_n-\bbbc(A_n-1)\right\},
$
where $A_n=\left(1+\frac{M_n}{2\zeta(3)}\right)^{1/3},$ $U_n$, $V_n$, $M_n$ are sample centered logarithmic moments
\begin{eqnarray*}
U_n&=\frac{1}{n}\sum_{j=1}^n\ln|Z_j|,
V_n=\frac{1}{n}\sum_{j=1}^n\left(\ln|Z_j|-U_n\right)^2,\\
M_n&=\frac{1}{n}\sum_{j=1}^n\left(\ln|Z_j|-U_n\right)^3,
\end{eqnarray*}
$\mathbf{I}(A)$ is the indicator of event $A$, $\bbbc=0.577\ldots$ is the Eulerian constant, and  $\zeta(3)$ is the Riemann function at point 3.

%
%


\begin{thebibliography}{}

\bibitem[Bening {\em et~al.}(2004)Bening, Korolev, Kolokoltsov, Uchaikin,
  Saenko, and Zolotarev]{Bening2004}
Bening, V.~E., Korolev, V.~Y., Kolokoltsov, V.~N., Uchaikin, V.~V., Saenko,
  V.~V., and Zolotarev, V.~M. (2004).
\newblock {Estimation of parameters of fractional stable distributions}.
\newblock {\em Journal of Mathematical Sciences\/}, {\bf 123}(1), 3722 -- 3732.

\bibitem[Bening {\em et~al.}(2006)Bening, Korolev, Sukhorukova, Gusarov,
  Saenko, Uchaikin, and Kolokoltsov]{Bening2006}
Bening, V.~E., Korolev, V.~Y., Sukhorukova, T.~A., Gusarov, G.~G., Saenko,
  V.~V., Uchaikin, V.~V., and Kolokoltsov, V.~N. (2006).
\newblock {Fractionally stable distributions}.
\newblock In V.~Y. Korolev and N.~N. Skvortsova, editors, {\em Stochastic
  Models of Structural Plasma Turbulence\/}, pages 175--244. Brill Academic
  Publishers, Utrecht.

\bibitem[Bunday(1984)Bunday]{Bunday1984}
Bunday, B. (1984).
\newblock {\em {Basic Optimization Methods}\/}.
\newblock Hodder Arnold.

\bibitem[Chambers {\em et~al.}(1976)Chambers, Mallows, and Stuck]{Chambers1976}
Chambers, J.~M., Mallows, C.~L., and Stuck, B.~W. (1976).
\newblock {A method for simulating stable random variables}.
\newblock {\em Journal of the American Statistical Association\/}, {\bf
  71}(354), 340--344.

\bibitem[Furusawa and Kaneko(2003)Furusawa and Kaneko]{Furusawa2003}
Furusawa, C. and Kaneko, K. (2003).
\newblock {Zipf's Law in Gene Expression}.
\newblock {\em Physical Review Letters\/}, {\bf 90}(8), 8--11.

\bibitem[Hoyle {\em et~al.}(2002)Hoyle, Rattray, Jupp, and Brass]{Hoyle2002}
Hoyle, D.~C., Rattray, M., Jupp, R., and Brass, A. (2002).
\newblock {Making sense of microarray data distributions.}
\newblock {\em Bioinformatics (Oxford, England)\/}, {\bf 18}(4), 576--84.

\bibitem[Kanter(1975)Kanter]{Kanter1975}
Kanter, M. (1975).
\newblock {Stable Densities Under Change of Scale and Total Variation
  Inequalities}.
\newblock {\em The Annals of Probability\/}, {\bf 3}(4), 697--707.

\bibitem[Kolokoltsov {\em et~al.}(2001)Kolokoltsov, Korolev, and
  Uchaikin]{Kolokoltsov2001}
Kolokoltsov, V.~N., Korolev, V.~Y., and Uchaikin, V.~V. (2001).
\newblock {Fractional Stable Distributions}.
\newblock {\em Journal of Mathematical Sciences\/}, {\bf 105}(6), 2569--2576.

\bibitem[Kotulski(1995)Kotulski]{Kotulski1995}
Kotulski, M. (1995).
\newblock {Asymptotic distributions of continuous-time random walks: A
  probabilistic approach}.
\newblock {\em Journal of Statistical Physics\/}, {\bf 81}(3-4), 777--792.

\bibitem[Kuznetsov {\em et~al.}(2002)Kuznetsov, Knott, and
  Bonner]{Kuznetsov2002}
Kuznetsov, V.~A., Knott, G.~D., and Bonner, R.~F. (2002).
\newblock {General statistics of stochastic process of gene expression in
  eukaryotic cells.}
\newblock {\em Genetics\/}, {\bf 161}(3), 1321--1332.

\bibitem[Liebovitch {\em et~al.}(2006)Liebovitch, Jirsa, and
  Shehadeh]{LIEBOVITCH2006}
Liebovitch, L.~S., Jirsa, V.~K., and Shehadeh, L.~A. (2006).
\newblock {Structure of genetic regulatory networks: evidence for scale free
  networks}.
\newblock In {\em Complexus Mundi - Emergent Patterns in Nature\/}, pages 1--8,
  Singapore. World Scientific Publishing Co. Pte. Ltd.

\bibitem[Lu and King(2009)Lu and King]{Lu2009}
Lu, C. and King, R.~D. (2009).
\newblock {An investigation into the population abundance distribution of
  mRNAs, proteins, and metabolites in biological systems.}
\newblock {\em Bioinformatics (Oxford, England)\/}, {\bf 25}(16), 2020--7.

\bibitem[Macneil and Walhout(2011)Macneil and Walhout]{Macneil2011}
Macneil, L.~T. and Walhout, A. J.~M. (2011).
\newblock {Gene regulatory networks and the role of robustness and
  stochasticity in the control of gene expression.}
\newblock {\em Genome research\/}, {\bf 21}(5), 645--57.

\bibitem[Metzler and Klafter(2000)Metzler and Klafter]{Metzler2000}
Metzler, R. and Klafter, J. (2000).
\newblock {The random walk's guide to anomalous diffusion: a fractional
  dynamics approach}.
\newblock {\em Physics Reports\/}, {\bf 339}(1), 1--77.

\bibitem[Saenko(2012)Saenko]{Saenko2012}
Saenko, V.~V. (2012).
\newblock {Maximum likelihood algorithm for approximation of local
  fluctuational fluxes at the plasma periphery by fractional stable
  distributions}.
\newblock {\em arxiv.org\/}, (arXiv:1209.2297 [physics.plasm-ph]).

\bibitem[Samko {\em et~al.}(1973)Samko, Kilbas, and Marichev]{Samko1973}
Samko, S.~G., Kilbas, A.~A., and Marichev, O.~I. (1973).
\newblock {\em {Fractional Integrals and Derivatives -Theory and
  Application}\/}.
\newblock Gordon and Breach, New York.

\bibitem[Uchaikin(2000)Uchaikin]{Uchaikin2000}
Uchaikin, V.~V. (2000).
\newblock {Montroll–Weiss problem, fractional equations, and stable
  distributions}.
\newblock {\em International Journal of Theoretical Physics\/}, {\bf 39}(8),
  2087--2105.

\bibitem[Uchaikin and Saenko(2002)Uchaikin and Saenko]{Uchaikin2002}
Uchaikin, V.~V. and Saenko, V.~V. (2002).
\newblock {Simulation of random vectors with isotropic fractional stable
  distributions and calculation of their probability density function}.
\newblock {\em J. Math. Sci.}, {\bf 112}(2), 4211 -- 4228.

\bibitem[Uchaikin and Zolotarev(1999)Uchaikin and Zolotarev]{Uchaikin1999}
Uchaikin, V.~V. and Zolotarev, V.~M. (1999).
\newblock {\em {Chance and stability Stable Distributions and their
  Applications}\/}.
\newblock VSP, Utrecht.

\bibitem[Ueda {\em et~al.}(2004)Ueda, Hayashi, Matsuyama, Yomo, Hashimoto, Kay,
  Hogenesch, and Iino]{Ueda2004a}
Ueda, H.~R., Hayashi, S., Matsuyama, S., Yomo, T., Hashimoto, S., Kay, S.~A.,
  Hogenesch, J.~B., and Iino, M. (2004).
\newblock {Universality and flexibility in gene expression from bacteria to
  human.}
\newblock {\em Proceedings of the National Academy of Sciences of the United
  States of America\/}, {\bf 101}(11), 3765--9.

\bibitem[Zolotarev(1986)Zolotarev]{Zolotarev1986}
Zolotarev, V.~M. (1986).
\newblock {\em {One-dimensional stable Distributions}\/}.
\newblock Amer. Mat. Soc., Providence, RI.

\end{thebibliography}

\end{document}